\documentclass[a4paper]{article}
\usepackage{latexsym}
\usepackage{amssymb}
\usepackage{amsmath}
\newtheorem{conj}{Conjecture}[section]
\newtheorem{lemma}[conj]{Lemma}
\newtheorem{coro}[conj]{Corollary}

\newtheorem{thm}{Theorem}

\newcommand{\qed}{\raisebox{-.8ex}{$\Box$}}
\newenvironment{bew}
{\noindent{\bf Proof:}}
{\hfill \qed\\}

\newcommand{\Aut}{{\rm Aut}}

\newcommand{\Syl}{{\rm Syl}}

\newcommand{\PSL}{{\rm PSL}}

\newcommand{\PSU}{{\rm PSU}}

\newcommand{\GL}{{\rm GL}}
\newcommand{\SL}{{\rm SL}}
\newcommand{\SU}{{\rm SU}}
\newcommand{\GU}{{\rm GU}}
\newcommand{\Sp}{{\rm Sp}}
\newcommand{\GF}{{\rm GF}}

\newcommand{\Stab} {{\rm Stab}}

\newcommand{\Alt}{~{\rm Alt}}

\newcommand{\diam}{~{\rm diam}}

\newcommand{\ZZ}{\mathbb{Z}}
\newcommand{\NN}{\mathbb{N}}
\newcommand{\QR}{{\rm I}\kern-5.0pt {\rm Q} \kern2pt}
\newcommand{\End}{{\rm End}}

\title{Commuting graphs of odd prime order elements in simple groups\thanks{This research is part of the project ``Transversals in Groups with an application to loops'' GZ: BA 2200/2-2 funded by
the DFG}}
\author{B. Baumeister, A. Stein}

\begin{document}

\maketitle

\begin{abstract}
We study the commuting graph on elements of odd prime order in finite simple groups.
The results are used in a forthcoming paper describing the structure of Bruck loops and Bol loops of exponent 2.
\end{abstract}

\section{Introduction}

Let $G$ be a group and $X$ a {\em normal subset of $G$}, that is for all $x \in X,g \in G$ we have 
$x^g \in X$. The {\em commuting graph} on $X$ is the undirected graph $\Gamma_{X,G}= \Gamma_X$ with vertex set $X$
such that two vertices $x$ and $y$, $x \ne y$, are on an edge if $[x,y]=1$. 
The commuting graph of a group is an  object which has been studied quite often to obtain strong results on the group
$G$. We give a short overview of some major work on or related to commuting graphs. For more details see the
references given below.

Bender noted in his paper on strongly 2-embedded subgroups, \cite{bender}, the equivalence between the existence of a strongly 2-embedded subgroup and the disconnectedness of the commuting graph of involutions.

At about the same time Fischer determined the groups generated by a class $X$
of $3$-transpositions by studying the commuting graph on $X$ [Fi].
Later Stellmacher classified those groups which are generated by a special class of 
elements of order $3$ again by examinating the related commuting graph [St].

In order to prove the uniqueness of the group of Lyons, Aschbacher and Segev showed that the commuting 
graph of $3$-subgroups generated by the $3$-central elements of a group of Lyons is simply connected [AS]. Notice also that a major breakthrough towards the famous Margulis-Platonov conjecture has been made by Segev by  using the commuting graph on the whole set $G$ for $G$ a non-trivial finite group [Se].

Finally Bates et all [BBPR] determined the diameter of the connected commuting graphs of a conjugacy class of involutions
of $G$ where $G$ is a Coxeter group and Perkins [Pe] did the same for the affine groups $\tilde{A}_n$, see also the 
related work [IJ2].
In [AAM] Abdollahi,  Akbari and Maimanithe considered the dual of the commuting graph
on $G\setminus{Z(G)}$. They conjectured that if these graphs are isomorphic for two non-abelian finite groups
then the groups have the same order. This conjecture has been checked for some simple groups in [IJ1].

In this paper given a finite simple group $G$, we consider the commuting graph $\Gamma_{\cal O}$ on the set 
${\cal O}$ of odd prime order elements of $G$  and some of its subgraphs. Our aim is to describe the connectivity of the graph. 

For some integer $n$ let $\pi(n)$ be the set of prime divisiors of $n$, for a group $G$ let $\pi(G):=\pi(|G|).$ 
Sometimes we consider the set of $\psi(n):=\pi(n)-\{ 2 \}$. 

For $G$ a group and $\rho$ a set of primes let ${\cal E}_\rho(G):=\{ x \in G| o(x) \in \rho \}$, (the set
of elements in $G$ of order a prime in $\rho$) and
for $X$ a normal subset of $G$ let $\Gamma_X$ be the commuting graph on $X$.
Notice, that $G$ acts by conjugation on $\Gamma_X$, inducing automorphisms. 
For $\rho \subseteq \pi(G)$ let $\Gamma_\rho= \Gamma_{{\cal E}_{\rho}(G)}$ and for an integer $n$ let $\Gamma_n=\Gamma_{\psi(n)}$.  Thus for $p$ a prime, $\Gamma_p$ is the commuting graph on the set of elements of order $p$ of $G$.

For $x \in X$ let ${\cal C}_x$ be the connected component of $\Gamma_X$ containing $x$ and $H_x \le G$ its stabilizer in $G$.
A connected component ${\cal C}_x$ is {\em big}, if  $G$ acts by conjugation on it, otherwise {\em small}. 
So, if ${\cal C}_x$ is big, then it contains the full conjugacy class $x^G$.

Let $q$ be a power of a prime $p$ and $r \ne p$ another prime. Set
\[ d_q(r):=\min \{ i \in \NN: r \mid q^i-1 \}.  \]
So $d_q(r)$ is the order of $q$ modulo $r$. 

In [S] as well as in [BS] we use the results of this paper to characterize the finite Bol loops of exponent $2$ as well as the finite Bruck loops.

Here in this paper we show

\begin{thm}
\label{theorem_big_component}
Let $G$ be a finite simple group. 
\begin{itemize}
\item[(a)] If $\Gamma_{\cal O}$ does not have a  big connected component, then $G$ is one of the following groups:
\begin{itemize}
\item $A_1(q), {}^2B_2(q), {}^2G_2(q)$ (for any $q$),
\item ${}^2A_2(q)$ for $q$ odd with $\frac{q+1}{(q+1,3)}$ a 2-power or
\item $M_{11},J_1,A_2(4)$
\end{itemize}
Conversely,  the mentioned groups do not have a big connected component. 
\item[(b)] If $\Gamma_{\cal O}$ has a big connected component, then this 
component contains an element $x$, $o(x)=p$, such that $\Gamma_p$ is connected.
In particular, by \ref{bender-case} and \ref{cyclic-sylows},  $G$ has no strongly $p$-embedded subgroup 
and the Sylow $p$-subgroups are not cyclic.
\end{itemize}
\end{thm}

This implies the following.

\begin{coro}
Let $G$ be a finite group and suppose that there is a big connected component $C$ of $\Gamma_{\cal O}$.
Suppose that the Sylow-$p$-subgroups are cyclic for all the primes $p$ such that there is an element of order $p$ in $C$.
Then $G$ is not simple.
\end{coro}

The uniqueness question of the big connected components is answered during the classification of small 
connected components below.

In the next theorem we give examples of big connected components which are a single conjugacy class.
We say that a {\em conjugacy class $x^G$ is connected} if the commuting graph on $x^G$ is connected.
In particular it follows that the commuting graph $\Gamma_{\cal O}$ of a group which has a connected conjugacy class $x^G$ with $x$ an element of
odd prime order does have a big connected component. In order to prove Theorem~\ref{theorem_big_component} for the alternating groups and the groups of Lie type in even characteristic, we give in the next theorem a list of elements $x$
for these groups, such that the commuting graph on $x^G$ is connected. This list is not complete, but contain all
alternating groups and groups of Lie type in even characteristic, which possess such an element $x$ of odd prime order. 

 \begin{thm}
\label{theorem_connected_conjugacy_class}
In the following table we list alternating groups and groups of Lie type in even characteristic $G$,
and subsets $\omega$ of  $\pi(G)$ such that the conjugacy class $x^G$ is connected for some element $x$ in $G$ of order $r$, 
  $r \in \omega$, with $E(C_G(x))/Z(E(C_G(x)))$ as given in the third column. In the first column of the table we list $G$, in the second $\omega$ is given and in the last further conditions which have to be satisfied. 
  
$$\begin{array}{|l|l|l|l|}
\hline	 	
G 				& \omega 			& E(C_G(x))/Z(E(C_G(x)))& \mbox{\rm conditions}\\
\hline
\Alt_n, n \ge 8 		& \{ 3 \} 			& \Alt_{n-3} 		& \\	
\hline	
A_2(q), q~\mbox{\rm even} 	&\pi(\frac{q-1}{(q-1,3)})	&  A_1(q)		& q>4 \\
 A_3(q), q~\mbox{\rm even}   			&\pi(q-1)  	& A_2(q) 		& q>2 \\
 A_n(q), q~\mbox{\rm even}, n\ge 4 		&\pi(q^2-1) 	& A_{n-2}(q) 	      	& \\
\hline
 {}^2A_2(q), q~\mbox{\rm even}			&\pi(\frac{q+1}{(q+1,3)})      	& A_1(q) 		& q>2 \\
 {}^2A_3(q), q~\mbox{\rm even} 			&\pi(q+1)	& {}^2A_2(q)		& \\
{}^2A_n(q), q~\mbox{\rm even}, n \ge 4		&\pi(q^2-1) 	&  {}^2A_{n-2}(q) 	& \\
\hline
 C_n(q), q~\mbox{\rm even}, n \ge 3		&\pi(q^2-1)	& C_{n-1}(q)		& \\
\hline
D_n(q), q~\mbox{\rm even}, n \ge 4		&\pi(q-1)	&  D_{n-1}(q)		& q>2 \\
 D_n(q), q~\mbox{\rm even}, n \ge 4		&\pi(q+1)	& {}^2D_{n-1}(q)	& \\
\hline
 {}^2D_n(q), q~\mbox{\rm even}, n \ge 4		&\pi(q-1)	& D_{n-1}(q)		& q>2 \\
 {}^2D_n(q), q~\mbox{\rm even}, n \ge 4		&\pi(q+1)	& {}^2D_{n-1}(q)	& \\
\hline
 G_2(q), q~\mbox{\rm even} 			&\pi(q^2-1)      & A_1(q)		& q>4 ~\mbox{\rm or}\\
						&		 &			& (q,r)=(4,3) \\ 	
\hline
 {}^3D_4(q), q~\mbox{\rm even}			&\pi(q^2-1) 	& A_1(q^3)		& q>2 \\
 {}^3D_4(2)					&\{3\}		& 3^{1+2}.2\Sigma_4 & \\
\hline
 {}^2F_4(q), q~\mbox{\rm even}			&\pi(q^2+1)      & {}^2B_2(q)		& q>2 \\
 {}^2F_4(2)'					&\{3\} 		 & 3^{1+2}.4 		&  \\
\hline
 F_4(q), q~\mbox{\rm even}				&\pi(q^2-1)	& C_3(q)		& \\	
\hline
 E_6(q), q~\mbox{\rm even}				&\pi(q^2-1)	& A_5(q)		& \\
\hline
 {}^2E_6(q), q~\mbox{\rm even} 			&\pi(q^2-1)	& {}^2A_5(q)		& \\
\hline
 E_7(q), q~\mbox{\rm even}				&\pi(q^2-1) 	& D_6(q)		& \\
\hline
 E_8(q), q~\mbox{\rm even}				&\pi(q^2-1)	& E_7(q)		& \\		
\hline
\end{array}$$
\end{thm}

\begin{thm}
\label{theorem_unique_big_connected_component}
Let $G$ be a finite simple group, such that $\Gamma_{\cal O}$ has a big connected component. 
Then either $G$ has a unique big connected component or $G=O'N$ and 
$\Gamma_{\{3,5\}}$ and $\Gamma_7$ are the two big connected components.
\end{thm}

In order to prove Theorem~\ref{theorem_unique_big_connected_component} we need the following information.

\begin{thm}
\label{theorem_small_connected_components}
Let $G$ be a finite simple group. Suppose that there is a big connected component and let $x$ be an element
 of $G$ of odd prime order $r$.
 If $x$ is not contained in a 
big connected component, then the conditions given in the table hold. Conversely if $x$ satisfies
these conditions, then $x$ is in a small connected component. 

$$\begin{array}{|l|l|}
\hline
G & r \\
\hline
M_{12} & r \in \{5,11\}  		\\
M_{22} & r \in \{5,7,11 \} \\
J_2    & r =7 							\\
M_{23} 	& r \in \{7,11,23 \} \\
HS     & r \in \{7,11\} 		\\
J_3    & r \in \{17,19\} \\
M_{24} & r \in \{11,23\} 		\\
 McL    & r \in \{7,11\} \\
He     & r =17 							\\
	Ru     & r\in \{ 7,13,29\} \\
Suz    & r\in \{ 11,13\} 		\\
	O'N    & r \in \{11,19,31\} \\
Co_3   & r \in \{ 11,23\} 	\\
	Co_2   & r \in \{ 7,11,23\} \\
Fi_{22}& r \in \{ 11,13 \} 	\\
 HN     & r \in \{11,19\} \\
Ly     & r \in \{ 31,37,67 \} \\
Th     & r \in \{ 19,31 \} \\
Fi_{23}& r \in \{ 11,17,23 \} \\
 Co_1   & r =23 \\
J_4    & r \in \{ 23,29,31,37,43\} \\
 Fi_{24}'& r \in \{ 17,23,29 \} \\
B      & r \in \{ 17,19,23,31,47 \} \\
	M      & r \in \{ 41,47,59,71 \} \\ 
\hline
\end{array}$$

$$\begin{array}{|l|l|l|}
\hline
G & ~\mbox{\rm condition on }~G & r\\
\hline
\Alt_n & n-t~\mbox{\rm a prime}~, t\in \{0,1,2\} & n-t \\
\hline
\end{array}$$

$$\begin{array}{|l|l|l|}
\hline
G & \mbox{\rm condition on }~G & d_q(r) \\
\hline
A_2(q)							&	\pi(\frac{q-1}{(q-1,3)}) \not\subseteq \{2 \}										  & 3\\
										& \pi(\frac{q-1}{(q-1,3)}) \subseteq \{ 2 \}, q~\mbox{\rm odd} 				  & 1,2,3 \\
A_3(q) 							& \pi(q+1) \subseteq \{ 2 \} 																			  & 4\\
          					& \pi(q-1) \subseteq \{ 2 \} 																			  & 3\\
A_n(q), n \ge 4 		&	q=3, n=4																												  & 4~(r=5)\\
                  	& n ~\mbox{\rm a prime},~ \pi(\frac{q-1}{(q-1,n+1)}) \subseteq \{ 2 \}  & n\\
                 		& n+1~\mbox{\rm a prime}																							  & n+1 \\
\hline
{}^2A_2(q)					& \pi(\frac{q+1}{(q+1,3)}) \not\subseteq \{ 2 \} 									  & 6 \\
{}^2A_3(q) 					& \pi(q-1) \subseteq \{ 2 \} 																			  & 4 \\
              			& \pi(q+1) \subseteq \{ 2 \} 																			  & 6 \\
{}^2A_n(q), n \ge 4	& q\in \{3,9 \}, n=4																								  & 4~ (r=5,41) \\
                    & n ~\mbox{\rm a prime}~, \pi(\frac{q+1}{(q+1,n+1)}) \subseteq \{ 2 \}  & 2n \\
                    & n+1   ~\mbox{\rm a prime}																						  & 2n+2 \\
\hline
B_n(q), n \ge 3, ~ q~\mbox{\rm odd} & \pi(n) \subseteq \{ 2\} 															& 2n \\
																& n~\mbox{\rm a prime},~ \pi(q-1) \subseteq \{ 2 \} 				& n \\
                    						& n~\mbox{\rm a prime},~ \pi(q+1) \subseteq \{ 2 \} 				& 2n \\
\hline
C_2(q)							& q \ne 2																													  & 4 \\
C_n(q), n \ge 3			& \pi(n) \subseteq \{ 2\} 																				  & 2n \\
                    & n~\mbox{\rm a prime},~ \pi(q-1) \subseteq \{ 2 \} 										& n \\
                    & n~\mbox{\rm a prime},~ \pi(q+1) \subseteq \{ 2 \} 										& 2n \\
\hline
D_n(q), n \ge 4			& n  ~\mbox{\rm a prime},~ \pi(q+1) \subseteq \{ 2 \} 									& n \\
                    & n-1~\mbox{\rm a prime},~ \pi(q-1) \subseteq \{ 2 \} 									& n-1 \\
                    & n-1~\mbox{\rm a prime},~ \pi(q+1) \subseteq \{ 2 \} 									& 2n-2\\
                    & \pi(n-1) \subseteq \{ 2 \}, \pi(q+1) \subseteq \{ 2 \} 					  & 2n-2\\
\hline
{}^2D_n(q),n \ge 4  & n~\mbox{\rm a prime},~ \pi(q+1) \subseteq \{ 2 \} 										& 2n \\
                    & \pi(n) \subseteq \{ 2 \}																				  & 2n \\
                    & n-1~\mbox{\rm a prime},~ q = 3																		& n-1, 2n-2 \\
                    & \pi(n-1) \subseteq \{ 2 \}, \pi(q-1) \subseteq \{ 2 \} 						& 2n-2 \\

\hline
G_2(q), q\ne 2 	& 3 \nmid q-1 														& 3 \\
              	& 3 \nmid q+1															& 6 \\
\hline              	
{}^3D_4(q)			& 																				& 12 \\
\hline
{}^2F_4(q)'			& 																				& 12 \\
								& q=2																			& 4(r=5)\\
\hline
F_4(q) 					& 																				& 8,12 \\
\hline
E_6(q)  				& 																				& 9 \\
        				& q \in \{ 3,7 \}																		& 8 (r=41,1201) \\
\hline        				
{}^2E_6(q) 			& 																				& 18 \\
        				& q\in \{ 2,3,5	\}																& 8 (r=17,41,313) \\
          			& q=2																			& 12 (r=13) \\
\hline          			
E_7(q)  				&  \pi(q+1) \subseteq \{ 2 \} 						& 14,18 \\
        				&  \pi(q-1) \subseteq \{ 2 \}							& 7,9 \\
\hline        				
E_8(q) 					&  																				& 15,24,30 \\
                &  5 \nmid q^2+1													& 20 \\
\hline                
\end{array}$$ 
\end{thm}

\begin{coro}\label{cyclic_small}
Let $G$ be a finite simple group and suppose $\Gamma_{\cal O}$ has a big connected component.
If $x$ is an element of $G$ of prime order $r$ which is in a small connected component, 
then $O^{2}(C_G(x))$ is abelian and the Sylow-$r$-subgroups of $G$ are either cyclic or $G \cong {}^2F_4(2)'$ and $r=5$.
\end{coro}

We wonder whether there is a proof, which does not use the full classification.
Notice, if $G =\PSL_2(8) \times Sz(8)$, then   $\Gamma_{\cal O}$ is connected and all the Sylow subgroups of odd order
are cyclic.

To prove our results we need studying the groups of Lie type closely.
We use the following sources about maximal subgroups of groups of Lie type: 
\cite{KL} for classical groups, \cite{LSS} and \cite{CLSS} for exceptional groups of Lie type.
Furthermore, the papers \cite{Coo}, \cite{K3D4} and \cite{Malle} were useful.

Our strategy to prove the theorems is as follows.
If $G$ is a finite simple group which is not listed in Theorem \ref{theorem_big_component} (a), then 
we prove Theorem~\ref{theorem_big_component} by either using the $p$-local subgroups of $G$, see 
Lemma~\ref{geometric_criterion}, or by producing
a connected conjugacy class in $G$ which shows Theorems~\ref{theorem_big_component} and \ref{theorem_connected_conjugacy_class} at the same time.  This is done in Section~4, first for the alternating, sporadic and then separately for
the groups of odd and even type, respectively.

Theorems~\ref{theorem_unique_big_connected_component} and \ref{theorem_small_connected_components}
are proven in Section~5.
Here our strategy is as follows: Let $C$ be a big connected component. Then we are able to show that
$C$ is the set of elements of $G$ of order $r$ with $r$ in $\rho$, for some $\rho \subseteq \pi(G)$.
The knowledge of centralizers and certain subgroups of $G$ then allows to describe $\rho$. The size of the centralizer of an element $x$ of $G$ then implies
whether $x$ is in  a small connected component or not, see Corollary~\ref{cyclic_small}.

In Section~2 we provide some facts from number theory and Section~3 contains general results about
commuting graphs and big connected components.

\section{Facts from number theory}

Let $q$ be a power of the prime $p$ and let $r \ne p$ be another prime. 
Recall $d_q(r) | r-1$ by Lagrange.
Let $n$ be an integer, $n \neq 0$. The famous theorem of K.Zsigmondy states

\begin{thm}\label{zsig}
There is either an odd prime $s$  with $d_q(s)=n$ 
or one of the following cases holds.

\begin{itemize}
\item[(a)] $q$ is a Mersenne prime, i.e. $q=p=2^m-1$ for some prime $m$ and $n=2$.
\item[(b)]  $q$ is a Fermat prime, i.e. $q=p=2^{2^m}+1$ for some integer $m$ or $q = 9$ and $n=1$.
\item[(c)]  $q=2$ and $n = 1$ or $n =6$ 
\end{itemize}
\end{thm}

Let $\Phi_n(x) \in \ZZ[x]$ be the $n$-th cyclotomic polynomial.
Then the following lemmata are consequences of Theorem~\ref{zsig}. 

\begin{lemma}
\label{Zsig}
Let $p$ be a prime and $n$ an integer. The following holds.
\begin{itemize}
\item[(a)] If $\Phi_n(p)$ a power of 2, then $n=1$ and $p$ is 2 or a Fermat prime or
$n=2$ and $p$ is a Mersenne prime.
\item[(b)]If $\Phi_n(p)$ is a power of 3, then $p=2$ and  $n \in \{1,2,6\}$.
\item[(c)]If $\Phi_n(p)$ a power of 3 times a power of 5, then $p=2$ and  $n \in \{1,2,4,6\}$.
\end{itemize}
\end{lemma}
\begin{bew}
If $n>2$ and $(p,n) \ne (2,6)$ by Zsigmondy's theorem 
there exists a prime $r$ dividing $\Phi_n(p)$, which does not divide $\Phi_m(p)$ for $m < n$.
Since $3$ divides $(p-1) p (p+1) = \Phi_1(p) p \Phi_2(p)$ we have $r>3$. 
So in the first two cases the question reduces to those primes $p$, 
for which $p-1$ (in case $n=1$)  or $p+1$ (in case $n=2$) is a 2-power or a 3-power. 
For the third case observe, that $n \mid r-1$, so $n \in \{1,2,4\}$ in this case and we have to determine
those primes $p$, for which one of $p-1,p+1$ or $p^2+1$ is a 3-power times a 5-power. Since in particular
$\Phi_n(p)$ is odd, $p=2$. The statement is immediate.
\end{bew}

\begin{lemma}
\label{qLemma}
Let $q$ be a prime power. The following holds.
\begin{enumerate}
\item[(a)] If $q-1$ is a 2-power, then $q=2$, $q=9$ or $q$ is a Fermat prime.
\item[(b)] If $q+1$ is a 2-power, then $q$ is a Mersenne prime. 
\item[(c)] If $q^2-1$ is a 2-power, then $q=3$.
\item[(d)] If $q^2-1$ is a 2-power times a 3-power, then $q \in \{2,3,5,7,17\}$. 
\item[(e)] If $q^2-1$ is a 3-power times a 5-power, then $q \in \{2,4\}$. 
\end{enumerate}
\end{lemma}
\begin{bew}
Let $q= p^e$. Remember the formulas 
\[ (p^e)^n -1 = \prod\limits_{d\mid en} \Phi_d(p)\] and 
\[ (p^e)^n+1 = \prod\limits_{\substack{d\mid 2en\\d\nmid en}} \Phi_d(p). \]
For $n=1$ we get $e \le 2$ in (i) and (ii) by  \ref{Zsig}.\\
For $n=2$ we get (iii) again by \ref{Zsig}.\\
Since 3 divides exactly one of $q-1$,  $q$, $q+1$, we get $q=2$ or $q$ a Mersenne or Fermat prime by (i) and (ii).\\
For Mersenne primes $p=2^r-1$ we have $p-1 = 2 (2^{r-1}-1)$, which is a 2-power times a 3-power
for $r\le2$ only by the formula mentioned and \ref{Zsig}.\\
For Fermat primes $p=2^m+1$ we can again use the formula on $p+1= 2(2^{m-1}+1)$ and \ref{Zsig}. 
Finally (v) is a consequence of the above product formula together with \ref{Zsig}.
\end{bew}

\section{Commuting graphs and big connected components}

We begin with some trivial but powerful observations. 
\begin{lemma}
\label{graphobservations}
Let $X$ be a normal subset of the group $G$ and $\Gamma_{X}$ the commuting graph on $X$.
\begin{itemize}
\item[(a)] $G$ acts by conjugation as a group of automorphisms on $\Gamma_X$.
\item[(b)] Let $g \in G$. Then the vertices $x^g$ and $x$ in $\Gamma_{X}$ are connected or equal
if and only if $g \in H_x$. 
\end{itemize}
\end{lemma}

The following lemma is helpful as it allows to switch from $G$ to a central extension of $G$.
 
\begin{lemma}
\label{center_quotient}
Let $X$ be a normal subset of the group $G$, $\Gamma_{X}$ the commuting graph on $X$ and $\overline{G}:=G/Z(G)$. 
If $x$ and $y$ are  elements in $X$ which are connected in $\Gamma_X$, then $\overline{x},\overline{y}$ are connected in $\Gamma_{\overline{X}}$.

\end{lemma}

\begin{lemma}
Suppose $C$ is a big connected component in $\Gamma_X$ which is a subset of  
$X:={\cal E}_{\rho}(G)$ for some $\rho \subseteq \pi(G)$.
If there is an element $x \in C$ of order $r$, then $C$ contains all the elements of order $r$.
\end{lemma}

\begin{bew}
Let $z \in X$ be of order $r$. 
We show, that $x$ and $z$ are connected in $\Gamma_X$. 
Let $R \in \Syl_r(G)$ with $z \in R$ and $g \in G$ with $y^g \in R$. 
Then $y^g$ and $z$ are connected via $Z(R)\ne 1$, as ${\cal E}_r(G) \subseteq X$. 
Therefore $(y,z^{g^{-1}})$,$(x,z^{g^{-1}})$ and $(x^g,z)$ are connected.
As ${\cal C}_x$ is big, $(x,x^g)$ are connected, so $(x,z)$ are connected. 
\end{bew}

\begin{coro}
\label{normalconnectedcomponents}
Let $\emptyset \ne X \subseteq {\cal O}$, such that $\Gamma_X$ is connected and such 
that $X^g=X$ for all $g \in G$. 
Then a subset $\rho \subseteq \pi(G)-\{2\}$ with $\{o(x): x \in X \} \subseteq \rho$ exists,
such that ${\cal E}_\rho(G)$ is the connected component in $\Gamma_{\cal O}$ containing $X$. 
In particular big connected components of $\Gamma_{\cal O}$ are subsets ${\cal E}_\rho(G)$. 
\end{coro}

Notice, that the subset $\rho$ for a big connected component $C$ of $\Gamma_\pi$ can be determined from the sizes of centralizers only,
once the order $r$ of a single element $x \in C$ is known. For this we simply define a graph on the set $\pi$ by connecting all primes $p_1$ and $p_2$,
such that $p_2$ divides the size of a centralizer of an element of order $p_1$. 
The connected component of the prime $r$ in this graph is the subset $\rho$ in question.

In order to use this method, we have to establish the existence of big connected components. 
A special case is the connectedness of $\Gamma_p$. Following Bender \cite{bender}, we show,
that connectedness of $\Gamma_p$ is equivalent to the fact that $G$ has no strongly $p$-embedded subgroup.
First we give criteria for the connectedness of $\Gamma_p$.

\begin{lemma}
\label{p-locals}
Let $G$ be a group with $O_p(G) \ne 1$. Then $\Gamma_p$ is connected.  
\end{lemma}

\begin{bew}
Choose $x \in \Omega_1(Z(O_p(G)))$. Now for $g \in G$ also $x^g \in \Omega_1(Z(O_p(G)))$,
so $[x,x^g]=1$ and $g \in H_x$. 
\end{bew}

\begin{lemma}
\label{geometric_criterion}
Suppose there exists a prime $p \in \pi(G)$ with $G = \langle N_G(Y): Y \le P, Y \ne 1 \rangle$ for some $P \in \Syl_p(G)$. 
Then $\Gamma_p$ is connected. 
\end{lemma}

\begin{bew}
Let $x \in \Omega_1(Z(P)), o(x)=p$. Then $P \le H_x$. 
For $1 \ne Y \le P$ we may choose $1 \ne y \in Y$ with $o(y) =p$. Then $N_G(Y) \le H_y$ by \ref{p-locals}.
As $H_x=H_y$, $H_x = \langle N_G(Y): Y \le P, Y \ne 1 \rangle = G$. Therefore all conjugates of $x$ in $G$ are connected,
so $\Gamma_p$ is connected.
\end{bew}

In some sporadic groups we need a generalization to include nonlocal subgroups $U$ with  $\Gamma_p(U)$ connected.\\

\begin{lemma}
\label{amalgam_criterion}
Suppose there exists a prime $p$ and subgroups $A,B \le G$, such that $G = \langle A,B \rangle$, $A \cap B$ contains elements of order $p$ and both $\Gamma_p(A)$ and $\Gamma_p(B)$ are connected.  
Then $\Gamma_p$ is connected.
\end{lemma}

\begin{bew}
Choose $x \in A \cap B, o(x)=p$. Consider $H_x$ in $\Gamma_p$. As $\Gamma_p(A)$ is connected, $A \le H_x$. As $\Gamma_p(B)$ is connected, $B \le H_x$.
Therefore $G=\langle A,B \rangle \le H_x$, so $\Gamma_x$ is connected.
\end{bew}

Recall that a subgroup $U \le G$ is {\bf strongly $p$-embedded}, if $U \neq G$, $p \in \pi(U)$ and $p \not\in \pi(U \cap U^g)$
for all $g \in G-U$, cf. \cite{bender}. The equivalence of (a) and (b) is already shown by Bender for $p = 2$ as 
well as essentially the equivalence of (b) and (c), see \cite{bender}.

\begin{lemma}
\label{bender-case}
Let $G$ be a finite group and $p \in \pi(G)$. The following statements are equivalent:
\begin{itemize}
\item[(a)] The graph $\Gamma_p$ is connected.
\item[(b)]  $G$ has no strongly $p$-embedded subgroup.
\item[(c)] For $P \in Syl_p(G)$: $G = \langle N_G(Y): 1 \ne Y \le P \rangle$.
\end{itemize}
\end{lemma}

\begin{bew}
Suppose $\Gamma_p$ is connected, but there exists a strongly $p$-embedded subgroup $U$. 
Let $x \in U, o(x)=p$. As $U$ is strongly $p$-embedded, $U$ is the stabilizer of a unique point in the action of $G$ on
the $U$-cosets and this is the unique fixed point of $x$. 
Therefore $C_G(x)$ fixes this unique point, so $C_G(y) \le U$ for every $y \in U$ of order $p$. 
This gives a contradiction to $\Gamma_p$ connected, as $G-U$ contains elements of order $p$. 

Suppose $U:=\langle N_G(Y): 1 \ne Y \le P \rangle \ne G$, but $G$ has no strongly $p$-embedded subgroup. 
Let $g \in G-U$ with $|U \cap U^g|_p$ maximal and $X \in \Syl_p(U \cap U^g)$. 
If $X=1$, then $U$ is strongly $p$-embedded, contrary to assumption. 

If $X \in \Syl_p(G)$, we find some $u \in U$ with $X^u = P$, so $U = \langle N_G(Y): 1 \ne Y \le X \rangle$. 
Likewise we find some $v \in U^g$ with $X^v=P^g$. Then also $U^g = \langle N_G(Y): 1 \ne Y \le X \rangle$, so $U=U^g$.
Then $g \in N_G(U)$. As $N_G(P) \le U$, $N_G(U)= U$ by Frattini, so $g \in U$, a contradiction.

So $1 < |X| < |G|_p$. Let $A,B \in \Syl_p(N_G(X))$ with $A \le U$ and $B \le U^g$. As $|A| > |X|$, $B \not\le U$.
We can choose a $Q \in \Syl_p(U)$ with $X \le Q$. There exists a $w \in U$ with $P^w=Q$, so $U = \langle N_G(Y): 1 \ne Y \le Q \rangle$. 
Then $N_G(X) \le U$ contradicts $B \not\le U$. 

Suppose $G= \langle N_G(Y): 1 \ne Y \le P \rangle$. Then $\Gamma_p$ is connected by \ref{geometric_criterion}. 
\end{bew}

\begin{coro}
\label{cyclic-sylows}
Let $G$ be a finite group and $p \in \pi(G)$. If $\Gamma_p$ is connected, then Sylow-$p$-subgroups of $G$ are noncyclic or $O_p(G) \ne 1$.
\end{coro}

\begin{bew}
As $\Gamma_p$ is connected, $G= \langle N_G(Y): 1 \ne Y \le P \rangle$. If Sylow-$p$-subgroups are cyclic, all those subgroups $N_G(Y)$ 
are contained in the subgroup $N_G(Y_1)$ for $Y_1=\Omega_1(P)$, so $O_p(G)$ contains $\Omega_1(P)$.
\end{bew}

In the groups of Lie type in even characteristic we wish to show connectedness of a conjugacy class.
the following lemma is a powerful tool.

\begin{lemma}
\label{UABg_criterion}
Let $\Gamma=\Gamma_X$ for $X=x^G, x \in {\cal O}$.  
Suppose  $U$ is a subgroup of $G$ such that $U = A B$ for two commuting subgroups $A$ and $B$ of $U$  and 
such that there is a $g \in G$ with $A^g \le B$.
Then $H_x \ge \langle U,g \rangle > U$.
\end{lemma}

\begin{bew}
Notice, that the commuting graph of $Y=X \cap U$ in $U$ is connected: $H_x(\Gamma_Y)$ contains $x$, so $B \le C_U(x)$,
so $A^g \le B$, so $x^g \in A^g$, so $A \le C_U(x^g)$, so $U$.
Furthermore, $g \in H_x=H_x(\Gamma)$, as $x$ and $x^g$ are connected in $U$. 
Therefore $\langle U,g \rangle \le H_x$. 
As $U \le N_G(A)$, but $g \not \in N_G(A)$, $\langle U,g \rangle > U$. 
\end{bew}

For $G$ a group of Lie type we show the existence of $U,A,B$ and $g$  by applying either the Curtis-Tits Theorem
or by using the Steinberg relations or by using the 
action of $G$ on its geometry or a natural module as is explained later. In some cases another criterion is useful:

\begin{lemma}
\label{abelian_criterion}
Let $G$ be a group and $x \in G$ an element of order $p$. 
If $G= \langle N_G(A) : A \le G, x \in A, A'=1 \rangle$, then $x^G$ is connected.  
\end{lemma}

\begin{bew}
Let $\Gamma=\Gamma_X$ for $X=x^G$. If $x \in A$ with $A'=1$, then $N_G(A) \le H_x$. So the condition implies $G \le H_x$,
so $\Gamma_X$ is connected. 
\end{bew}

We end this section with a criterion for the nonexistence of big connected
components. 

\begin{lemma}
\label{nonabelian_centralizer}
Let $X$ be a normal subset of $G$, and $C$ a big connected component of $\Gamma_X$. 
Then either some  $x \in C$ exists such that $C_G(x)$ is not abelian or $\langle C \rangle \le F(G)$. 
\end{lemma}

\begin{bew}
Suppose $C_G(x)$ is abelian for every $x \in C$. Let $x,y,z \in C$ with
$[x,y]=1=[y,z]$. As $C_G(y)$ is abelian and $x,z \in C_G(y)$, $[x,z]=1$. 
As $C$ is a connected component, any two elements of $C$ commute, so
$A:=\langle C \rangle$ is abelian. As $C$ is a big connected component, $A$
is $G$-invariant, so $A \le F(G)$.
\end{bew}

Notice, that groups with abelian centralizers were considered already by L.Weisner \cite{weisner} and 
M.Suzuki \cite{suz-CA-groups}.
We wonder, whether it is possible to classify those finite simple groups without big connected component in 
$\Gamma_{\cal O}$ without using the classification.

\section{Proofs of Theorems~\ref{theorem_big_component} and \ref{theorem_connected_conjugacy_class}}

In this section we show that if $G$ is a simple group not listed in Theorem~\ref{theorem_big_component} (a),
then $\Gamma_{\cal O}$ has at least one big connected component. At the same time we show that there is
a prime $p$ such that $\Gamma_p$ is connected.

For groups of Lie type in even characteristic the strategy is to establish 
Theorem~\ref{theorem_connected_conjugacy_class}. This, then produces big connected components by 
 \ref{normalconnectedcomponents}. 

\subsection{Alternating groups}

\begin{lemma}
\label{pcycles}
Let $G \cong \Alt_n$ and $x \in G$ of odd prime order $p$. 
\begin{itemize}
\item[(1)] $O_p(C_G(x))$ contains $p$-cycles.
\item[(2)] If $x$ is a $p$-cycle, then:
\begin{itemize}
\item[(a)] If $p+p < n$, then the commuting graph on $x^G$ is connected.
\item[(b)] $F^\ast(C_G(x)) \cong \langle x \rangle \times A_{n-p}$, unless $n-p=4$.
\item[(c)] If $p$ is not a Fermat prime, then $|N_G(\langle x \rangle):C_G(x)|$ is divisible by some odd prime $r$ dividing $p-1$.
\item[(d)] If $p+3\le n$, then $C_G(x)$ contains a 3-cycle.
\end{itemize}
\end{itemize}
\end{lemma}
\begin{bew} 
The centralizer of an element of order $p$ acts on the fixed points and permutes the cycles of lenght $p$.
This gives (1),(2b) and (2d). For (2c) we observe, that in $\Sigma_n$ all powers of $x$ are conjugate, as they have the same cycle structure.
Remains (2a): For a $p$-cycle $x$ let $M(x)\subseteq \{ 1,...,n\}$ be the orbit of length $p$. 
Now, if for $p$-cycles $x,y$: $|M(x) \cap M(y)|=p-1$, then $x,y$ are connected in the commuting graph:
Since $|M(x)\cup M(y)| = p+1 \le n-p$, some $p$-cycle $z$ exists with $M(x) \cap M(z) = \emptyset = M(y) \cap M(z)$,
so $[x,z]=1=[y,z]$. But now, given any two $p$-cycles $x,y$,
we can find $p$-cycles $z_i$ with: $z_0:=x$, $z_k=y$ and $|M(z_i) \cap M(z_{i+1})|=p-1$ for $0 \le i<k$. Therefore
the commuting graph on $x^G$ is connected.
\end{bew}

\begin{lemma}
Theorem \ref{theorem_connected_conjugacy_class} holds for $G$ an alternating group. This implies
that Theorem \ref{theorem_big_component} holds for $G$, as well.
\end{lemma}

\begin{bew}
By \ref{pcycles} (2)(a), the conjugacy class of 3-cycles is connected for $n \ge 7$.
\end{bew}

\subsection{Sporadic groups}

\begin{lemma}
Theorem \ref{theorem_big_component} holds for $G$ a sporadic group.
\end{lemma}

\begin{bew}
We use the informations from \cite{ATLAS}. 
Notice, that $M_{11}$ and $J_1$ have no big connected component, as visible from the centralizer sizes.

We use \ref{geometric_criterion} to establish the connectedness of $\Gamma_p$ in the following cases $(G,p)$:

$(M_{12},3)$, $(J_2,5)$, $(J_3,3)$, $(McL,3)$, $(He,7)$, $(Ru,5)$, $(Suz,3)$, $(O'N,3)$, $(Co_3,3)$, $(Fi_{23},3)$, $(Co_1,3)$, $(Fi_{24}',3)$, $(B,3)$, $(M,3)$.

We use \ref{amalgam_criterion} with subgroups $A,B$ in the following cases $(G,o(x),A,B)$:\\
$(M_{22}, 3, \Alt_7,\Alt_7)$, $(M_{23},3,M_{22},\Sigma_8)$, $(HS,3,M_{22},\Alt_8)$, $(M_{24},3,3\Sigma_6,\Alt_8)$, 
$(O'N,7,A_2(7):2,A_2(7):2)$, $(Co_2,3,McL,3^{1+4}:2^{1+4}.\Sigma_5)$, $(J_4,3,6M_{22},M_{24})$.
\end{bew}

\subsection{Groups of Lie type in odd characteristic}

\begin{lemma}
\label{odd_pcc}
Theorem \ref{theorem_big_component} holds for $G$ a group of Lie type in characteristic
$p>2$ other than $A_1(q), {}^2A_2(q)$ or ${}^2G_2(q)$. 
\end{lemma}

\begin{bew}
In this case $G$ is generated by its $p$-locals. Then \ref{geometric_criterion} 
shows, that $\Gamma_p$ is connected. 
\end{bew}

\begin{lemma}
\label{PSL2odd}
Theorem \ref{theorem_big_component} holds for $G \cong A_1(q)$, $q$ odd.
\end{lemma}

\begin{bew}
From Dixon's theorem on subgroups of $\PSL_2(q)$ we conclude the centralizer
sizes, so no big component exists. Also \ref{nonabelian_centralizer} implies
the statement. 
\end{bew}

\begin{lemma}
Theorem \ref{theorem_big_component} holds for $G \cong {}^2A_2(q)$, $q$ odd.
\end{lemma}

\begin{bew}
Let $q=p^e$ with $p$ a prime. Given an element $x \in G$, $o(x)=r$ an odd prime, we have either $r = p$,
$r \mid q+1$, $r \mid q-1$ or $r \mid \frac{q^2-q+1}{(3,q+1)}$. 
A Sylow-$p$-subgroup is strongly $p$-embedded, so $\Gamma_p$ is not connected.
By \ref{nonabelian_centralizer} we look for nonabelian centralizers. The only nonabelian centralizer of 
a semisimple element is a subgroup of type $\frac{q+1}{(q+1,3)} \circ \SL_2(q).2$. 

Suppose that $\frac{q+1}{(q+1,3)}$ is a 2-power.
Then all semisimple elements have an abelian centralizer. The nonabelian
centralizers from elements of ${\cal O}$ come from elements of order $p$.
But if $x,y$ are elements in $\Gamma_O$, $o(x)=p$, $[x,y]=1$  and $y$
semisimple, then $C_G(y) \le C_G(x)$, so $\diam({\cal C}_x)=1$. 
Let $N:=\langle {\cal C}_x \rangle \le C_G(x)$. As ${\cal C}_x$ should be a big connected
component, $N \trianglelefteq G$, a contradiction to the simplicity of $G$. 

Suppose now, that $\frac{q+1}{(q+1,3)}$ is not 2-power. 
Let $x \in G$ be an element with centralizer isomorphic to $\frac{q+1}{(q+1,3)} \circ \SL_2(q).2$.
We claim, that the conjugacy class $X:=x^G$ is connected. 
We can find $x$ in an abelian subgroup $A$ of size $\frac{(q+1)^2}{(q+1,3)}$. 
Consider $H_x$ in $\Gamma_X$. Clearly $C_G(x)$ is contained in it. 
But also $N_G(A) \le H_x$, as for $g \in N_G(A)$: $[x,x^g]=1$. As $N_G(A)/A \cong \Sigma_3$ and $C_G(x)$
is a maximal subgroup of $G$ not containing $N_G(A)$, $G = \langle  C_G(x),N_G(A) \rangle \le H_x$, 
so $x^G$ is connected. 
\end{bew}

\begin{lemma}
Theorem \ref{theorem_big_component} holds for $G \cong {}^2G_2(q)'$, $q$ odd.
\end{lemma}

\begin{bew}
The case $q=3$ is treated as $\PSL_2(8)$ in the next section.

We use the list of maximal subgroups in \cite{K2G2}. 
In particular as centralizers of semisimple elements are reductive,
centralizers of elements of odd prime order in $G$ are abelian $3'$-groups.
So centralizers of elements of order 3 are $\{2,3\}$-groups. 

By \ref{bender-case} therefore $\Gamma_3$ is not connected.
Then by \ref{nonabelian_centralizer}, $G$ has no big connected component. 
\end{bew}

\subsection{Groups of Lie type in even characteristic}

We first establish Theorem \ref{theorem_connected_conjugacy_class}, as then
Theorem \ref{theorem_big_component} is a consequence of \ref{normalconnectedcomponents}.

\begin{lemma}
Theorem \ref{theorem_big_component} holds for $G \cong A_1(q)$ and $G \cong {}^2B_2(q)$, $q>2$ even.
\end{lemma}

\begin{bew}
Use Dixon's Theorem for $\PSL_2(q)$ and \cite{suzuki} in case of $Sz(q)$ for the list of maximal subgroups.
Then \ref{nonabelian_centralizer} shows, that $G$ has no big connected component.
\end{bew}

\begin{lemma}
\label{PSL3connected}
\label{PSL4connected}
Let $G \cong A_n(q)$ with $n =2$ or $n=3$ and $q$ even.
Then Theorem \ref{theorem_connected_conjugacy_class} and
Theorem \ref{theorem_big_component} hold for $G$.
\end{lemma}

\begin{bew}
If $q=2$, the group $A_2(2)\cong A_1(7)$ has no big connected subgroup by \ref{PSL2odd}. As $A_3(2) \cong \Alt_8$,
we have a big connected component by \ref{pcycles}.

The group $A_2(4)$ has no big connected component, as visible from the centralizer sizes in \cite{ATLAS}.

So let $q\ge 4$ ($q>4$ for $n=2$) and $r$ a prime divisor of $q-1$. 
Let $a,b,c \in \GF(q)$ with $1 \ne a$, $a^r=1$, $b=\frac{1}{a^2}$ and $c=\frac{1}{a^3}$. 

For $n=2$ let $x_1$ the image of ${\rm Diag}(a,a,b)$ in $G$ and $x_2$ the image of ${\rm Diag}(b,a,a)$.
Then $[x_1,x_2]=1$, $x_1,x_2$ are conjugate in $G$ and $\langle x_1,x_2 \rangle \cong \ZZ_r \times \ZZ_r$.
We calculate, that with $q_1:=\frac{q-1}{(q-1,3)}$ we have $C_G(x_1) \cong \ZZ_{q_1} \times \PSL_2(q)$, $C_G(x_2) \cong \ZZ_{q_1} \times \PSL_2(q)$ and  
$N_G(\langle x_1,x_2 \rangle) \cong (\ZZ_{q_1} \times \ZZ_{q-1}):\Sigma_3$. 

Let $\Phi=\{ r_1,r_2,r_1+r_2, -r_1,-r_2,-r_1-r_2 \}$ be a root system  of type $A_2$
such that $\SL_3(q) = \langle x_r(t)~:~r \in \Phi, t \in GF(q)$ and such that 

$C_G(x_1)$ contains the image of $\langle X_{r_1},X_{-r_1} \rangle $,

$C_G(x_2)$ contains the image of $\langle X_{r_2},X_{-r_2} \rangle$ and

and $N_G(\langle x_1,x_2 \rangle)$ contains the image of the subgroup $N$, see \cite{Car}. 

The commutator relations imply that $\langle C_G(x_1),C_G(x_2) \rangle=G $. 
As $N$ is transitively on roots, it follows that $\langle C_G(x_1),N_G(\langle x_1,x_2 \rangle) \rangle \ge \langle C_G(x_1),C_G(x_2) \rangle=G$. 
Hence, $x_1^G$ is connected by \ref{abelian_criterion}.

For $n=3$ let $y_1$ the image of ${\rm Diag}(a,a,a,c)$ in $G$, $y_2$ the image of ${\rm Diag}(c,a,a,a)$ in $G$
and $y_3$ the image of ${\rm Diag}(a,c,a,a)$ in $G$.

We calculate that $[y_1,y_2]=1=[y_1,y_3]=[y_2,y_3]$, the $y_1,y_2,y_3$ are conjugate in $G$ and 
$\langle y_1,y_2,y_3 \rangle \cong \ZZ_r \times \ZZ_r \times \ZZ_r$. Moreover for $d=(q-1,3)$ we get
$C_G(y_1) \cong \ZZ_{q-1} . \PSL_3(q).\ZZ_d$, $C_G(y_2) \cong \ZZ_{q-1} . \PSL_3(q).\ZZ_d$,
$C_G(y_3) \cong \ZZ_{q-1} . \PSL_3(q).\ZZ_d$ and $N_G(\langle y_1,y_2,y_3 \rangle) \cong (\ZZ_{q-1} \times \ZZ_{q-1} \times \ZZ_{q-1}):\Sigma_4$.

Again we take a root system $\Phi$ of type $A_3$ for $\SL_4(q)$ with fundamental system $\Pi=\{ r_1,r_2,r_3\}$, such that
$C_G(y_1)$ contains the image of $\langle X_{r_2},X_{-r_2}, X_{r_3},X_{-r_3} \rangle$
and $N_G(\langle y_1,y_2,y_3 \rangle) = N$. 

We get $\langle C_G(y_1),N_G(\langle y_1,y_2,y_3 \rangle) \rangle \ge \langle X_{r_1},X_{-r_1}, X_{r_2},X_{-r_2}, X_{r_3},X_{-r_3} \rangle = G$,
so by \ref{abelian_criterion} $y_1^G$ is connected.
\end{bew}

\begin{lemma}
\label{PSLconnected}
Let $G \cong A_n(q)$ for $n \ge 4$ and $q$ even.

Then Theorem \ref{theorem_connected_conjugacy_class} and
Theorem \ref{theorem_big_component} hold for $G$. 
\end{lemma}

\begin{bew}
We show the statement in $\hat{G}=\SL_n(q)$, which stays valid in $G= \hat{G}/Z(\hat{G})$ by \ref{center_quotient}.
Let $\Phi$ be a root system of $\hat{G}$ with fundamental root set $\Pi=\{ r_1,r_2,...,r_n\}$ with usual numbering of roots,
as described in \cite{Car}. 

Set $U= \langle X_r,X_{-r},  r \in \Pi-\{ r_2 \} \rangle$, $A= \langle X_{r_1},X_{-r_1} \rangle$ and $B= \langle X_r,X_{-r},  r \in \Pi-\{ r_1,r_2 \}$.
Let $x \in A$ be some element of order $r$ for $r$ a prime divisor of $q^2-1$. 

Choose $g \in N$, such that $g$ acts on $\Pi$ as transposition $(r_1,r_n)$. 
By \ref{UABg_criterion} then $\langle U,g \rangle \le H_x$. 
But $\langle U,g \rangle$ contains $\langle X_r,X_{-r},  r \in \Pi \rangle = \hat{G}$.
So $x^G$ is connected in $\hat{G}$. 
\end{bew}

In the unitary case let $F=\GF(q^2)$, $\alpha : a \mapsto a^q \in \Aut(F)$. 

Let $n>2$ be some integer and extended $\alpha$ to $\GL_n(F)$. For $g \in \GL_n(F)$ let $g^T$ be the transpose of $g$. 
Then \[ \GU_n(q)= \{ g \in \GL_n(F): g g^{ \alpha T} =1 \}\] and \[ \SU_n(q)=\{ g \in \GU_n(q): \det(g)=1 \}.\] 
A diagonal matrix ${\rm Diag}(a_1,a_2,...,a_n)$ is contained in $\SU_n(q)$, iff $a_i^{q+1}=1$ for all $i$ and $\prod_{i=1}^n a_i=1$. 
Recall, that $\GU_n(q)$ preserves the unitary form $(u,v)= \sum\limits_{i=1}^n u_i v_i^\alpha$ for $u=(u_1,u_2,...,u_n) \in F^n$, $v=(v_1,v_2,...,v_n) \in F^n$.

For $u,v \in F^n$ we have $u \perp v$ iff $(u,v)=0$. For $U \le F^n$ let $U^\perp =\{ v \in F^n: u \perp v ~\mbox{for all}~ u \in U \} \le F^n$.

Let $e_i = (a_1,a_2,...,a_n) \in F^n$ be the standard base with $a_j=\delta_{i,j}$.
For $u \in F^n$ let $N(u)=(u,u) \in \GF(q)$ and  observe $N(\lambda u)= \lambda^{q+1} N(u)$. 
By Hilbert 90 , the map $F \to \GF(q): a \mapsto a^{q+1}$ is surjective, so for every $u \in F^n$ some $\lambda \in F$ exists
with $N(\lambda u)=1$. 
Recall, that $\SU_n(q)$ acts transitively on orthonormal bases.

\begin{lemma}
\label{PSU3connected}
\label{PSU4connected}
Let $G \cong {}^2A_2(q)$ for $n=2$ or $n=3$ and $q$ even.
Then Theorem \ref{theorem_connected_conjugacy_class} and
Theorem \ref{theorem_big_component} hold for $G$.  
\end{lemma}

\begin{bew}
If $q=2$, the group ${}^2A_2(2)$ is soluble. The group ${}^2A_3(2)$ is isomorphic to $B_2(3)$, so by \ref{odd_pcc} has a big connected component.

So let $q>2$ and $r$ a prime divisor of $q+1$. There exist $a,b,c \in F$ with $a^r=1$, $b=\frac{1}{a^2}$ and $c=\frac{1}{a^3}$.
If $q=8$ and $n=2$ choose $1 \ne a$, $a^9=1 \ne a^3$ and $b=\frac{1}{a^2}$.

We do our calculations in $\SU_n(q)$ and use \ref{center_quotient} for the proof of the statement.

For $n=2$ let $x_1={\rm Diag}(a,a,b)$ and  $x_2= {\rm Diag}(b,a,a)$ in $G$. We will show, that $x_1^{\hat{G}}$ is connected in $\SU_3(q)$.
We can calculate, that $o(x_1)=r$, $x_1,x_2 \in \hat{G}=\SU_3(q)$, $[x_1,x_2]=1$, $x_1,x_2$ are conjugate in $\SU_3(q)$, $A_1:=C_{\hat{G}}(x_1) \cong \ZZ_{q+1} \times \PSL_2(q)$
and $A_2:=N_{\hat{G}}(\langle x_1,x_2 \rangle) \cong (\ZZ_{q+1} \times \ZZ_{q+1}):\Sigma_3$. By \ref{abelian_criterion}, $H_{x_1} \ge \langle A_1,A_2 \rangle=:G_0$.
Furthermore \[ A_1= \Stab_{\hat{G}}(\langle e_3 \rangle) = \Stab_{\hat{G}}(\langle e_1,e_2 \rangle) \]  and 
\[ B_1:=C_{\hat{G}}(x_2)= \Stab_{\hat{G}}(\langle e_1 \rangle) = \Stab_{\hat{G}}(\langle e_2,e_3 \rangle) \le G_0. \] 
We can now show, that $G_0 = \hat{G}$:

Let $g \in \hat{G}$ and $v_i=e_i^g$, $i=1,2,3$. Let $u \in \langle v_1 \rangle^\perp \cap \langle e_1 \rangle^\perp, N(u)=1$. 
We can find some $g_1 \in B_1$ with $u^{g_1}=e_3$. As $v_1 \perp u$, $v_1^{g_1} \perp u^{g_1}=e_3$, so $v_1^{g_1} \in \langle e_1,e_2 \rangle$. 

We can find some $g_2 \in A_1$ with $v_1^{g_1 g_2} = e_1$.  
As $v_2,v_3 \in \langle v_1 \rangle^\perp$, $v_2^{g_1 g_2}, v_3^{g_1 g_2} \in \langle e_2,e_3 \rangle$. So there exists some $g_3 \in B_1$ with $v_i^{g_1 g_2 g_3} \in \langle e_i \rangle$,
so $g g_1 g_2 g_3 $ is a diagonal matrix. As $A_1 \le G_0$ contains all diagonal matrices of $\hat{G}$ and $g_1 g_2 g_3 \in G_0$, $g \in G_0$, so $\hat{G}=G_0$. 
Therefore $x_1^{\hat{G}}$ is connected.

For $n=3$ let $y_1={\rm Diag}(a,a,a,b)$, $y_2={\rm Diag}(b,a,a,a)$ and $y_3={\rm Diag}(a,b,a,a)$.
We will show, that $y_1^{\hat{G}}$ is connected in $\SU_4(q)$.
We can calculate, that $o(y_i)=r$, $y_i \in \SU_4(q)$, $[y_i,y_j]=1$, the $y_i$ are conjugate in $G$, 
$A_1:=C_{\hat {G}}(y_1) \cong \GU_3(q)$ and $A_2:=N_{\hat{G}}(\langle y_1,y_2,y_3 \rangle) \cong (\ZZ_{q+1} \times \ZZ_{q+1} \times \ZZ_{q+1}):\Sigma_4$.
Using the same method as in case $n=2$, we can show, that $\langle A_1,A_2 \rangle = \hat{G}$. By \ref{abelian_criterion}, $y_1^{\hat{G}}$ is connected.
\end{bew}

\begin{lemma}
\label{PSUconnected}
Let $G \cong {}^2A_n(q)$ for $n \ge 4$ and $q$ even.
Then Theorem \ref{theorem_connected_conjugacy_class} and Theorem \ref{theorem_big_component} hold for $G$. 
\end{lemma}

\begin{bew}
By \cite{KL} there exists a maximal subgroup $U$ of type $U_2(q) \perp U_{n-2}(q)$. Let $A,B$ be the subgroups
of $U$ with $A \cong \SL_2(q)$ and $B \cong \SU_{n-2}(q)$. 

As $\SU_n(q)$ acts transitively on nondegenerated 2-subspaces of its natural
module, there exists some $g \in G$, such that $A^g \subseteq B$. 

Let $r$ be any prime divisor of $q^2-1$ and $x \in A$ be some element of order $r$.
Using the list of maximal subgroups in \cite{KL} we conclude, that $C_G(x) = C_U(x)$. 
By \ref{UABg_criterion}, $H_x \ge \langle U,g \rangle$. As $U$ was maximal, $H_x =G$, so $x^G$ is connected.
\end{bew}

\begin{lemma}
\label{Sp4connected}
Let $G \cong C_2(q)$ for $q>2$, $q$ even. 
Then Theorem \ref{theorem_big_component} hold for $G$.  
\end{lemma}

\begin{bew}
Let $r \in \pi(q^2-1)$. We show, that $\Gamma_r$ is connected.
There exist two classes of maximal subgroups $M_1$,$ M_2$ of type $(\PSL_2(q) \times \PSL_2(q)).2$,
which are interchanged by a graph automorphism. 

We can choose $M_1$ to be of type $(\Sp_2(q) \perp \Sp_2(q)):2$, the normalizer of a 2-space decomposition and $M_2$ to be of type $O_4^+(q)$.

Notice, that these two subgroups contain Sylow-subgroups for all primes dividing $q^2-1$.
By Sylow's Theorem we may choose $M_1,M_2$ with a common Sylow-$r$-subgroup.
Let $x \in M_1 \cap M_2$ be of order $r$. Notice, that $\Gamma_r(M_1)$ and $\Gamma_r(M_2)$ are connected.
By \ref{amalgam_criterion}, $\Gamma_r(\langle M_1,M_2 \rangle)$ is connected. As $M_1,M_2$ are maximal subgroups,
$\Gamma_r$ is connected.
\end{bew}

\begin{lemma}
\label{connected_symplectic_CC}
Let $G \cong C_n(q)$ for $n \ge 3$, $q$ even. 
Then Theorem \ref{theorem_connected_conjugacy_class} and Theorem \ref{theorem_big_component} hold for $G$. 
\end{lemma}

\begin{bew}
By \cite{KL} there exists a maximal subgroup $U$ of type $\Sp_2(q) \perp \Sp_{2n-2}(q)$. Let $A,B$ be the normal subgroups
of $U$ with $A \cong \SL_2(q)$ and $B \cong \Sp_{2n-2}(q)$. 

As $\Sp_{2n}(q)$ is transitive on nondegenerate 2-spaces, there exists some $g \in G$,
such that $A^g \subseteq B$. 

Let $r$ be any prime divisor of $q^2-1$ and $x \in A$ be some element of order $r$.
By \ref{UABg_criterion}, $H_x \ge \langle U ,g \rangle >U$. As $U$ was a maximal subgroup, 
$\Gamma_X$ is connected for $X=x^G$. 
\end{bew}

\begin{lemma}
\label{Omega_connected}
Let $G \cong D_n(q)$ or ${}^2D_n(q)$ for $n \ge 4$, $q$ even.
Then Theorem \ref{theorem_connected_conjugacy_class} and Theorem \ref{theorem_big_component} hold for $G$. 
\end{lemma}

\begin{bew}
Let $\varepsilon \in \{+,-\}$.
By \cite{KL} there exist maximal subgroups $U^\varepsilon_+$ of type $O_2^+(q) \perp O^{\varepsilon}_{2n-2}(q)$ 
and $U^\varepsilon_{-}$ of type $O_2^-(q) \perp O^{-\varepsilon}_{2n-2}(q)$ in $\Omega^\varepsilon(q)$, 
provided $q>2$ in case $U^\varepsilon_+$. For $q=2$ we exclude the cases $U^\varepsilon_+$, as then $q-1=1$.

Let $A^\varepsilon_+ \cong \ZZ_{q-1}, B^\varepsilon_{+} \cong \Omega_{2n-2}^\varepsilon(q)$ be normal subgroups
of $U^\varepsilon_+$ and 
$A^\varepsilon_- \cong \ZZ_{q+1}, B^\varepsilon_{-} \cong \Omega_{2n-2}^{-\varepsilon}(q)$ be normal subgroups
of $U^\varepsilon_-$. 
From the action of $G$ on its natural module we conclude, that in any case some $g^{\varepsilon_1}_{\varepsilon_2}$
exist, such that $(U^{\varepsilon_1}_{\varepsilon_2},A^{\varepsilon_1}_{\varepsilon_2},B^{\varepsilon_1}_{\varepsilon_2},\End
g^{\varepsilon_1}_{\varepsilon_2})$
satisfy the conditions of \ref{UABg_criterion} for $x$ any element of order $r$, $r \in \pi(q-(\varepsilon_2 1))$. 
As $U^{\varepsilon_1}_{\varepsilon_2}$ was a maximal subgroup of $G$, $\Gamma_X$ is connected for $X=x^G$.
\end{bew}

\begin{lemma}
\label{G2connected}
Let $G \cong G_2(q)$ with $q$ even. 
Then Theorem \ref{theorem_connected_conjugacy_class} and Theorem \ref{theorem_big_component} hold for $G$. 
\end{lemma}

\begin{bew}
As $G_2(2)'\cong {}^2A_2(3)$, Theorem \ref{theorem_big_component} holds for $G_2(2)'$.

Let $q \ge 4$. We use the list of maximal subgroups in \cite{Coo}.
Let $\varepsilon \in \{ +,-\}$ with $r$ a divisor of $q- \varepsilon$.
There exist two classes of subgroups of type $(q-\varepsilon) \times \PSL_2(q)$
in a maximal subgroup of type $\PSL_2(q) \times \PSL_2(q)$. Let $C_1,C_2$ be representatives of the two classes
and $x_1 \in Z(C_1)$, $x_2 \in Z(C_2)$ with $o(x_1)=r=o(x_2)$.

Notice, that there is only one class of maximal subgroups $M$ isomorphic to $A^\varepsilon_2(q).2 \cong \SL^\varepsilon_3(q).2$ for each $\varepsilon$. 
We can choose $i \in \{1,2\}$, such that $M$ does not contain a conjugate of $C_i$, as $M$ contains a unique class
of such subgroups. 

Now $H_{x_i}$ contains $C_i$, but also a subgroup $N$ of shape $(q-\varepsilon)^2:D_{12} \le M$. 
So $H_{x_i} \ge \langle C_i,N \rangle$. Using the list of maximal subgroups of $G$, we
see for $q>4$, that $\langle C_i,N \rangle  \ge G$, as $C_i$ is not in a conjugate of $M$. Therefore $\Gamma_X$ for
$X=x_i^G$ is connected. 

Notice, that our selection of $C_i$ also forces $C_G(x_i) \cong (q-\varepsilon) \times \PSL_2(q)$, even if $r=3$.

In case $q=4$, due to the $J_2$-maximal subgroup, we used computer calculations.
We calculated in MAGMA, using the 6-dimensional representation of $G$ over $\GF(4)$,
that $G$ has connected conjugacy classes of elements of order 3 with the given centralizer structure.
There is no connected conjugacy class of elements of order 5, though $\Gamma_5$ is connected.
\end{bew}

\begin{lemma}
\label{3D4connected}
Let $G \cong {}^3D_4(q)$  with $q$ even. 
Then Theorem \ref{theorem_connected_conjugacy_class} and Theorem \ref{theorem_big_component} hold for $G$. 
\end{lemma}

\begin{bew}
For $q=2$ we use the list of maximal subgroups in \cite{ATLAS}. 
By \ref{geometric_criterion}, $\Gamma_3$ and $\Gamma_7$ are connected. As $G$ has 3 $3$-local maximal subgroups,
but only two classes of elements of order 3, $G$ has a connected conjugacy class of elements of order 3 
by \ref{abelian_criterion}. However it is class 3B, which is not the class used in case of $q>2$.

So let $q>2$ and $r$ any prime divisor of $q^2-1$. Let $\varepsilon \in \{-1,+1\}$
with $r$ a divisor of $q-\varepsilon$. 

We use the list of semisimple centralizers and maximal subgroups in \cite{K3D4}.

We can choose some $x \in G$, $o(x)=r$ with $C_G(x) \ge \ZZ_{q-\varepsilon} \times \PSL_2(q^3)$.
From the list of maximal subgroups we conclude, that $C_G(x)$ is contained in maximal parabolic
or a subgroup of type $\PSL_2(q) \times \PSL_2(q^3)$, as no other subgroup contains a $\PSL_2(q^3)$.
As centralizers of semisimple elements are reductive, $C_G(x) \cong \ZZ_{q-\varepsilon} \times \PSL_2(q^3)$. 
We claim that $\Gamma_X$ for $X=x^G$ is connected. Notice, that $G$ contains a torus normalizers $N$ of
type $\ZZ_{q^3-\varepsilon} \times \ZZ_{q-\varepsilon}. D_{12}$ in a subgroup $M$
of type $(q^2+\varepsilon q +1).A^\varepsilon_2(q).f_\varepsilon.2$ with $f_\varepsilon=(3,q-\varepsilon)$. 
As $M$ contains a Sylow-$r$-subgroup, we may assume $x \in M$ and $x \in N$. Notice, that $N$ is neither
contained in a maximal parabolic subgroup nor a subgroup of type $\PSL_2(q)\times \PSL_2(q^3)$. 

Therefore $\langle C_G(x),N \rangle = G$. By \ref{abelian_criterion} then $H_x \ge G$, so $\Gamma_{X}$ is connected.
\end{bew}

\begin{lemma}
\label{2F4connected}
Let $G \cong {}^2F_4(q)'$. 
Then Theorem \ref{theorem_connected_conjugacy_class} and Theorem \ref{theorem_big_component} hold for $G$. 
\end{lemma}

\begin{bew}
If $q=2$ we use this list of maximal subgroups in \cite{ATLAS}. 
By \ref{amalgam_criterion} we have $\Gamma_3$ connected. Notice, that $\Gamma_5$ is not connected,
as a Sylow-5-subgroup is normal in the centralizer of a 5-element.

For $q>2$ we use the list of maximal and maximal local subgroups in \cite{Malle}. 
Notice, that $5 \mid q^2+1$ in this case.

We can factorize $q^2+1 = (q-\sqrt{2q}+1)(q+\sqrt{2q}+1)$. Let $\varepsilon \in \{+,-\}$,
such that $r$ is a divisor of $q+\varepsilon \sqrt{2q}+1$ and let $x \in G$ be an element of order $r$
with $C_G(x) \cong \ZZ_{q+\varepsilon \sqrt{2q}+1} \times {}^2B_2(q)$. Such an element exists
in a maximal subgroup $M_1$ of type $({}^2B_2(q) \times {}^2B_2(q)).2$. Notice, that the outer involution
interchanges the components, as ${}^2B_2(q)$ has no outer automorphism of order 2. 
This gives $M_1 \le H_x$.

But there exists a subgroup $N$ of type $(\ZZ_{q+\varepsilon \sqrt{2q}+1} \times \ZZ_{q+\varepsilon \sqrt{2q}+1}).[96]$,
which is maximal for $q>8$ or $r>5$, while contained in ${}^2F_4(2)$ for $q=8$ and $r=5$. 
As $N \not\le M_1$, $H_x \ge G$, $\Gamma_X$ is connected for $X=x^G$. 
\end{bew}

\begin{lemma}
\label{F4connected}
Let $G \cong F_4(q)$ for $q$ even. 
Then Theorem \ref{theorem_connected_conjugacy_class} and Theorem \ref{theorem_big_component} hold for $G$. 
\end{lemma}

\begin{bew}
Let $r$ be a prime divisor of $q^2-1$.
By \cite{LSS}, $G$ has two classes of maximal subgroups $M_1,M_2$ isomorphic to $\Sp_8(q) \cong C_4(q)$.

By \ref{connected_symplectic_CC}, each $M_i$ has a connected conjugacy class for a prime $r \mid q^2-1$. 

We may choose $x \in M_1$ of order $r$ with $C_G(x)= C_{M_1}(x) \cong (q-\varepsilon) \times \Sp_6(q)$ for
for some $\varepsilon \in \{+,-\}$. 
The fact, that $C_G(x) = C_{M_1}(x)$ comes from the list of maximal subgroups, which contain a centralizer, 
see the main theorem of \cite{CLSS}.\\
Then $x$ is contained in a torus $T$ of type $(q-\varepsilon)^4$, with $W(F_4)$, the full Weyl group, acting on it.  
As this torus normalizer is not contained in $\Sp_8(q)$ (but in $\Omega^+_8(q).\Sigma_3$), we have $H_x =G$:
$H_x$ contains $M_1$ as seen in \ref{connected_symplectic_CC} and $N_G(T)$, but $\langle M_1,N_G(T) \rangle = G$,
as $M_1$ is a maximal subgroup not containing $N_G(T)$. Therefore the commuting graph on $x^G$ is connected.
\end{bew}

\begin{lemma}
\label{E678connected}
Let $G \cong E_6(q), {}^2E_6(q)$, $E_7(q)$ or $E_8(q)$ for $q$ even.
Then Theorem \ref{theorem_connected_conjugacy_class} and Theorem \ref{theorem_big_component} hold for $G$. 
\end{lemma}

\begin{bew}
By \cite{LSS} there are maximal subgroups $U$ with normal subgroups $A \cong \PSL_2(q)$
and $B \cong \PSL_6(q), \PSU_6(q), \Omega^+_{12}(q)$ resp. $E_7(q)$, 
such that a $g \in G$ exists with $A^g \subseteq B$. The existence of $g$
and these subgroups can also be seen from the Steinberg relations. 

Let $r$ be a prime divisor of $q^2-1$ and $x \in U$ be some element of order $r$. 
We can conclude from the main result of \cite{CLSS}, that $C_U(x)=C_G(x)$.
By \ref{UABg_criterion} the graph $\Gamma_X$ for $X=x^G$ is connected.
\end{bew}

\section{Proof of Theorem \ref{theorem_small_connected_components}}
We consider only those groups, which have a big connected component. 
Groups without big connected component were determined in Theorem \ref{theorem_big_component}.

We classify the small connected components and show uniqueness of the big connected components (if possible).
To do this, we use \ref{normalconnectedcomponents}: \\
We start with the prime(s) mentioned in the proof of Theorem \ref{theorem_big_component}.
We chose several big centralizers to show connectedness of a large subset of $\Gamma_{\cal O}$. 

We use knowledge on centralizers to show, that the remaining primes give elements, 
which are not connected to the big connected component. 

\subsection{alternating and sporadic groups}

\begin{lemma}
Theorem \ref{theorem_small_connected_components} holds for $G$ an alternating group.
\end{lemma}

\begin{bew}
This is a consequence of \ref{pcycles}. 
\end{bew}

\begin{lemma}
Theorem \ref{theorem_small_connected_components} holds for $G$ a sporadic group.
\end{lemma}

\begin{bew}
By Theorem \ref{theorem_big_component}, we can exclude $M_{11}$ and $J_1$.

By the centralizer sizes in \cite{ATLAS}, all primes listed as a small connected component give a unique small connected component. 

It remains to show, that the big connected component(s) contains all other primes. 
In the list below we give the set of primes $\pi({\cal C})$ of the orders of elements in the big  connected component 
together with elements $x$ whose centralizer size shows that the elements of $G$ of order $r$, $r$
in $\pi({\cal C})$, form indeed a connected component of $\Gamma_{\cal O}$.

This also shows, that the big connected component is unique, apart from the case $G=O'N$. 

$$
\begin{array}{|l|l|l|}
\hline
\mbox{Group} & \pi({\cal C}) & x \\
\hline
M_{12}	& \{ 3 \} 							&	3A \\
M_{22}  & \{ 3 \} 							& 3A \\
J_2			& \{ 3,5 \} 						& 3A \\
M_{23}  & \{ 3,5 \} 						& 3A \\
HS			& \{ 3,5 \}							& 3A \\
J_3			& \{ 3,5 \}							& 3A \\ 
M_{24}  & \{ 3,5,7 \} 					& 3A,3B \\
McL			& \{ 3,5 \}							& 3A \\
He			&	\{ 3,5,7 \}						& 3A \\
Ru			& \{ 3,5 \}							& 3A \\
Suz			& \{ 3,5,7 \}						& 3A \\
O'N			& \{ 3,5 \}							& 3A \\
				& \{ 7 \}								& 7A \\
Co_3		& \{ 3,5,7 \}						& 3A,3C \\
Co_2 		& \{ 3,5 \}							& 3A \\
Fi_{22} & \{ 3,5,7 \} 					& 3A \\
HN			& \{ 3,5,7 \}						& 5A \\
Ly			& \{ 3,5,7,11 \} 				& 3A \\
Th			& \{ 3,5,7,13 \}				& 3A,3C \\
Fi_{23} & \{ 3,5,7,13 \} 				& 3A \\
Co_1		& \{ 3,5,7,11,13 \}			& 3A \\
J_4			& \{ 3,5,7,11 \}				& 3A \\
Fi_{24}'& \{ 3,5,7,11,13 \}			& 3A,3B \\
B				& \{ 5,5,7,11,13 \}			& 3A \\
M				& \{ 3,5,7,11,13,17,19,23,29,31 \} & 3A,3C \\	
\hline	
\end{array}$$

\end{bew}

\subsection{Groups of Lie type}

We exclude the groups listed in Theorem \ref{theorem_big_component}. 

\begin{lemma}
\label{psl3uni}
Let $G \cong A_2(q)$. Then Theorem \ref{theorem_small_connected_components} holds.
\end{lemma}

\begin{bew}
Notice, that the torus of size $\frac{q^2+q+1}{(q-1,3)}$ is always self centralizing, so gives always a small connected component.\

If $q$ is even, $q>4$ by Theorem \ref{theorem_big_component}.
Then by Theorem \ref{theorem_connected_conjugacy_class} there is a connected conjugacy class with centralizer $\frac{q-1}{(q-1,3)} A_1(q)$. 
Therefore the big connected component is unique and consists of all elements of order $r$ with $r \in \pi(q^2-1)$. 

If $q$ is odd, by \ref{odd_pcc},  ${\cal E}_p(G)$ is connected. 
Centralizers of semisimple elements are either tori or of type $\frac{q-1}{(q-1,3)}\cdot L_2(q).2$. 
If $\frac{q-1}{(q-1,3)}$  is a 2-power, the big connected component contains no semisimple elements.
Else we find some element $x \in G$, $o(x)=r$ for some odd prime $r \ne p, r \mid q-1$
such that $C_G(x)$ contains a component isomorphic to $\SL_2(q)$. This shows $\Gamma_{(q-1)(q)(q+1)}$ connected. 
\end{bew}

\begin{lemma}
\label{psl4uni}
Let $G \cong A_3(q)$. Then Theorem \ref{theorem_small_connected_components} holds.
\end{lemma}

\begin{bew}
For $q=2$ we use $A_3(2) \cong \Alt_8$. 

If $q$ is even, $q \ne 2$, so by Theorem \ref{theorem_connected_conjugacy_class} there is a connected conjugacy class with centralizer $(q-1) A_2(q)$.
This shows $\Gamma_{(q^2-1)(q^3-1)}$ connected. There exists a subgroup $\ZZ_\frac{q^4-1}{q-1}$ from the $\GL_2(q^2) \le \GL_4(q)$. 
This subgroup contains elements of order $q+1$ in its center, so $\Gamma_{\cal O}$ is connected.

If $q$ is odd, $\Gamma_p$ is connected by \ref{odd_pcc}. There exists a subgroup of type $L_2(q) \oplus L_2(q)$, which shows,
that $\Gamma_{q (q^2-1)}$ is connected.

If $q-1$ is not a 2-power, we find a subgroup of type $L_1(q) \oplus L_3(q)$ with center of odd order. 
If $q-1$ is a 2-power, elements of order $r$ with $d_q(r)=3$ are contained in small connected components, as visible in $G \cong \Omega_6^+(q)$. 

If $q+1$ is not a 2-power, we find a subgroup of type $\GL_2(q^2) \le \GL_4(q)$ with center of odd order. 
If $q+1$ is a 2-power, elements of order $r$ with $d_q(r)=4$ are contained in small connected components, as visible in $G \cong \PSL_4(q)$. 
\end{bew}

\begin{lemma}
\label{psl5+uni}
Let $G \cong A_n(q)$ for $n \ge 4$. Then Theorem \ref{theorem_small_connected_components} holds.
\end{lemma}

\begin{bew}
Let $q$ odd.
There exists a subgroup of type $L_2(q) \oplus L_{n-1}(q)$. By \ref{odd_pcc}, $\Gamma_p$ is connected, so also
$\Gamma_{|\PSL_{n-1}(q)|}$ is connected.

If $q$ is even we get $\Gamma_{|\PSL_{n-1}(q)|}$ connected by Theorem \ref{theorem_connected_conjugacy_class}.
So remains to check the primes $r$ with $d_q(r)=n$ and $d_q(r)=n+1$. 
 
Suppose $d_q(r)=n+1$. 
If $n+1$ is a prime, a torus of size $\frac{q^{n+1}-1}{(q-1)(q^{n+1}-1,n+1)}$ is self centralizing and gives a small connected component.

If $n+1=a \cdot b$ with $a \ne 1 \ne b$, there exists a subgroup $M_1$ of type $L_{(n+1)/b}(q^b)$ in class ${\cal C}_3$. 

If $Z(F^\ast(M_1))$ contains elements of odd prime order, $\Gamma_{\cal O}(M_1)$ is connected.
As $F^\ast(M_1)$ contains a section isomorphic to $\PSL_2(q)$ and a torus of type $q^{n+1}-1$, 
elements of order $r$ with $d_q(r)=n+1$ are then contained in the big connected component.

By Proposition 4.3.6 of \cite{KL}, this subgroup is local with a cyclic normal subgroup of size $\frac{(q-1,b)(q^b-1)}{(q-1)(q-1,(n+1))}$. 
By Zsigmondy, some odd prime $t \mid q^b-1$ exists with $d_q(t)=b$, unless $b=2$ and $q$ is a Mersenne prime.
If $n$ is a 2-power, then $n\ge 8$ and there exists a subgroup $M_2 \le M_1$ of type $L_{{n+1}/4}(q^4)$.

Now $Z(F^\ast(M_2))$ has elements of odd order, as there exists a Zsigmondy-prime $t$ with $d_q(t)=4$. 
As $F^\ast(M_2)$ contains a torus of type $q^{n+1}-1$ and a $PSL_2(q)$-section, again elements of order $r$ with $d_q(r)=n+1$ are contained in the big connected component.

Suppose now $d_q(r)=n$. There exists a subgroup $M_3$ of type $L_1(q) \oplus L_n(q)$. 
By Proposition 4.1.4 of \cite{KL}, $Z(F^\ast(M_4))$ contains elements of odd order $s$ with $s \mid q-1$, 
if $\frac{q-1}{(q-1,n+1)}$ is not a 2-power. 
In that case $F^\ast(M_3)$ contains a torus of type $q^n-1$.

If $Z(F^\ast(M_3))$ contains no elements of odd prime order, $F^\ast(M_3)$ contains a component of type $L_n(q)$. 
The connected components of the commuting graph for $F^\ast(M_3)$ were determined by induction.
We have to distinguish the case $n=4$ where we use \ref{psl4uni} and $n>4$. 

For $n=4$ we have to care for small connected components of $F^\ast(M_3)$ containing elements of prime order $r$ with $d_q(r)=4$.
By \ref{psl4uni}, such connected components exist only for odd $q$ and then $q+1$ is a 2-power. As $F(M_3)=1$, $\frac{q-1}{(q-1,5)}$ is a 2-power too.
As neither $q-1$ nor $q+1$ is divisible by 3, $p=3$ and $q$ is a 3-power. If $q>3$, then $q+1$ has a Zsigmondy divisor bigger than 5, a contradiction.
The case $q=3$ arises.

For $n>4$ we have use induction. Again we have to care for small connected components of $F^\ast(M_3)$, which contain elements of prime order $r$
with $d_q(r)=n$. This forces $n$ to be a prime.

Notice, that from the action of $\SL_{n+1}(q)$ on its natural module, it is obvious,
that the listed cases all occure.
\end{bew}

\begin{lemma}
\label{psu3uni}
Let $G \cong {}^2A_2(q)$. Then Theorem \ref{theorem_small_connected_components} holds.
\end{lemma}

\begin{bew}
Notice, that for all $q$ a torus of size $\frac{q^2-q+1}{(q+1,3)}$ is self centralizing.

If $q$ is even, by Theorem \ref{theorem_connected_conjugacy_class} there is a connected conjugacy class $y^G$, $o(y)=r$ for $r$ some prime divisor of $q+1$. 
By construction of $y$, ${\cal C}_y$ contains ${\cal E}_{\psi(q^2-1)}(G)$, the statement holds for $q$ even.

So let $q$ odd and suppose there exist a big connected component, so $\frac{q+1}{(q+1,3)}$ is not a 2-power. Then a semisimple element of odd order $y$ exists,
such that $C_G(y) \cong \frac{q+1}{(q+1,3)} \circ \SL_2(q).2$. Therefore again a big connected component exists, containing ${\cal E}_\psi(q^2-1)(G)$. 
\end{bew}

\begin{lemma}
\label{psu4uni}
Let $G \cong {}^2A_3(q)$. Then Theorem \ref{theorem_small_connected_components} holds.
\end{lemma}

\begin{bew}
Consider first $q$ even.
By Theorem \ref{theorem_connected_conjugacy_class}, there is a connected conjugacy class $y^G$, $o(y)=r$ for $r$ some prime divisor of $q+1$. 
By construction of $y$, ${\cal C}_y$ contains ${\cal E}_{\psi(|\PSU_3(q)|)}(G)$.  
There exists a subgroup $\ZZ_\frac{q^4-1}{q+1}$ in a Levi complement of a parabolic subgroup of type $q^4: \GL_2(q^2)$.
This subgroup contains elements of order $s$ for $s$ some prime divisor of $q-1$, if $q>2$, so $\Gamma_{\cal O}$ is connected.
The case $q=2$ gives a small connected component.

So $q$ is odd. There exists a subgroup $M_1$ of type $U_2(q) \perp U_2(q)$. 
By \ref{odd_pcc}, a big connected component containing all elements of odd prime order $s$ with $s$ a divisor of $|\PSL_{n-2}(q)|$ exists. 
So remain the cases $d_q(r) \in \{ 4,6 \}$. 

In case of $d_q(r)=6$, let $M_2$ be a subgroup of type $U_1(q) \perp U_3(q)$. The
structure of $M_2$ is described by Proposition 4.1.4 of \cite{KL}. In particular $Z(F^\ast(M_2))$
contains elements of odd order, if $\frac{q+1}{(q+1,4)}$ is not a 2-power.  Notice, that $M_2$ contains a torus of type $q^3+1$.
Also for $x \in M_2$, $o(x)=r$ with $d_q(r)=6$, $C_{M_2}(x) = C_G(x)$ as visible from the action on the 6-dimensional $\GF(q)$-module
of $\Omega^-_6(q)$. Therefore we get a small connected component for $q+1$ a 2-power. 

If $d_q(r)=4$, let $M_3$ be a maximal subgroup of type $\GL_2(q^2)$ in class ${\cal C}_2$. 
The structure of $M_3$ is described by Proposition 4.2.4 of \cite{KL}. 
In particular $Z(F^\ast(M_3))$ has size $\frac{(q-1)(q+1,2)}{(q+1,4)}$, so contains elements of odd prime order, if $q-1$ is not a 2-power.
Notice that $M_3$ contains a torus of type $q^2+1$. Also for $x \in M_2$, $o(x)=r$ with $d_q(r)=4$, $C_{M_2}(x) = C_G(x)$ as visible from the action on the 4-dimensional $\GF(q^2)$-module
of $\SU_4(q)$. Therefore we get a small connected component for $q-1$ a 2-power. 
\end{bew}

\begin{lemma}
\label{psu5+uni}
Let $G \cong {}^2A_n(q)$ for $n \ge 4$. Then Theorem \ref{theorem_small_connected_components} holds.
\end{lemma}

\begin{bew}
If $q$ is even, by Theorem \ref{theorem_connected_conjugacy_class}, there is a connected conjugacy class $y^G$, $o(y)=r$ for $r$ some prime divisor of $q^2-1$.

If $q$ is odd, by \ref{odd_pcc}, a big connected component exists, containing all elements of order $p$. 

There exists a subgroup of type $U_2(q) \perp U_{n-1}(q)$. Therefore a big connected component exists,
which contains all elements of odd prime order $s$ with $s$ a divisor of $|\PSU_{n-1}(q)|$. 
So $r \mid (q^{n+1}-(-1)^{n+1}) (q^n-(-1)^n)$.

Suppose $n+1$ even and $r \mid q^{n+1}-1$. 
There exists a torus of type $q^{n+1}-1$ in a subgroup $M_1$ of type $\GL_{\frac{n+1}{2}}(q^2).2$ in class ${\cal C}_2$. 
If $\frac{n+1}{2}$ is even, then $\frac{n+1}{2} \ge 4$. Let $t$ be some Zsigmondy prime with $d_q(t)=4$. 

If $\frac{n+1}{2}$ is odd and $(q,n+1)\ne (2,6)$, let $t$ be some Zsigmondy prime with $d_q(t)=\frac{n+1}{2}$. 
If $(q,n+1)=(2,6)$ let $t=3$. Now the torus of type $q^{n+1}-1$ (and size $\frac{q^{n+1}-1}{(q+1)(q^{n+1}-1,n+1)}$ 
contains elements of order $t$, but $t \mid |\SU_{n/2}(q)|$, so $x$ is in the big connected component.

Suppose $n+1$ odd, but not a prime and $r \mid q^{n+1}+1$. Let $n+1=a \cdot b$ with $a \ne 1 \ne b$ and $b$ a prime.
There exists a subgroup $M_2$ of type $U_{\frac{n+1}{b}}(q^b)$ in class ${\cal C}_3$. 

By Proposition 4.3.6 of \cite{KL}, this subgroup is local with a cyclic normal subgroup of size $\frac{(q+1,b)(q^b+1)}{(q+1)(q+1,n+1)}$. 

By Zsigmondy, some odd prime $t \mid q^b-1$ exists with $d_q(t)=b$.
So $Z(F^\ast(M_2))$ contains elements of odd prime order, while $F^\ast(M_2)$ contains a $\PSL_2(q)$-section and 
a torus of type $q^{n+1}+1$. Therefore $x$ is contained in the big connected component.
If $n+1$ is a prime, a torus of type $q^{n+1}+1$ is self centralizing, so gives a small connected component. 

Suppose now $r \mid q^n-(-1)^n$. There exists a subgroup $M_3$ of type $U_1(q) \oplus U_n(q)$. 
By Proposition 4.1.4 of \cite{KL}, $Z(F^\ast(M_3))$ contains elements of odd order $s$ with $s \mid q+1$, 
if $\frac{q+1}{(q+1,n+1)}$ is not a 2-power. 
In that case $F^\ast(M_3)$ contains a torus of type $q^n-(-1)^n$, so $x$ is contained in the big connected component.

If $Z(F^\ast(M_3))$ contains no elements of odd prime order, $F^\ast(M_3)$ contains a component of type $U_n(q)$. 
We use the knowledge about the commuting graph of that component, but have to distinguish the case $n=4$, where we use \ref{psu4uni} and the case $n\ge 5$.

If $n=4$, we have to care for small connected components of $F^\ast(M_3)$, which contain elements of prime order $r$
with $d_q(r)=4$. By \ref{psu4uni} this makes $q-1$ a 2-power. As $\frac{q+1}{(q+1,5)}$ is a 2-power too,
$p=3$. For $q>9$, $q+1$ has a Zsigmondy divisor bigger than 5. So $q=3$ or $q=9$. 

If $n>5$, we have to care for small connected components of $F^\ast(M_3)$, which contain elements of prime order $r$
with $d_q(r)=2n$. This force $n$ to be a prime. The corresponding torus in $F^\ast(M_3)$ is self centralizing. 

Notice, that the listed small connected components arise, as visible from the action of $\SU_{n+1}(q)$ on its natural module.
\end{bew}

\begin{lemma}
\label{psp4uni}
Let $G \cong C_2(q)$. Then Theorem \ref{theorem_small_connected_components} holds.
\end{lemma}

\begin{bew}
By Theorem \ref{theorem_big_component} we have $q>2$.\\
Notice, that self centralizing tori of size $\frac{q^2+1}{(q-1,2)}$ exist, which give small connected components.
If $q$ is even, by \ref{Sp4connected}, a big component containing ${\cal E}_{\pi(q^2-1)}(G)$ exists.

If $q$ is odd, by \ref{odd_pcc} there exists a big connected component, containing all elements of order $p$. 
There exists a subgroup of type $\Sp_2(q) \perp \Sp_2(q)$. Therefore, if $r \mid (q-1)q(q+1)$, then $x$ is in the big connected 
component. 
\end{bew}

\begin{lemma}
\label{psp6+uni}
Let $G \cong C_n(q)$ for $n \ge 3$. Then Theorem \ref{theorem_small_connected_components} holds.
\end{lemma}

\begin{bew} 
If $q$ is odd, by \ref{odd_pcc} there exists a big connected component, containing all elements of order $p$. 
If $q$ is even, by Theorem \ref{theorem_connected_conjugacy_class} there is a big connected component
containing all elements of prime order $s$ for $s$ a divisor of $q^2-1$.
There exists a subgroup of type $\Sp_2(q) \perp \Sp_{2n-2}(q)$. 
Therefore, if $r \mid |\Sp_{2n-2}(q)|$, then $x$ is in the big connected component.

So $r \mid (q^n-1)(q^n+1)$. 
If $n$ is even, then $r \mid q^n+1$, else $\Sp_n(q)$ contains elements of order $r$.
Let $n=a \cdot b$ with $a$ a 2-power and $b$ odd. If $b=1$, we have a self centralizing torus of size $\frac{q^a+1}{(q-1,2)}$, which gives a
small connected component.  

So $b>1$. There exists a subgroup $M_1$ of type $\Sp_{2b}(q^a)$. This subgroup contains a subgroup $M_2$ of type $\GL_b(q^a)$, which contains a torus of type $q^n-1$, and 
$M_3$ of type $\GU_b(q^a)$, which contains a torus of type $q^n+1$. The structure of $M_2$ is described
by Proposition 4.2.5, while those of $M_3$ is described by 4.3.7 for $q$ odd and 4.3.18 for $q$ even. 

We see, that $Z(F^\ast(M_2))$ contains no elements of odd order, iff $q-1$ is a 2-power and $a=1$.
Furthermore $Z(F^\ast(M_3))$ contains no elements of odd order, iff $q+1$ is 2-power and $a=1$. 
Both subgroups contain a $\PSL_2(q)$-section. 

If $a>1$, then $x$ is in the big connected component. 
If $a=1$, $b$ a prime and $q-1$ is a 2-power, we have a small connected component to a torus of size $q^b-1$.
If $a=1$, $b$ a prime and $q+1$ a 2-power, we have another small connected component to a torus of size $q^b+1$.
The existence of these small connected components is visible from the natural module of $\Sp_{2n}(q)$.

Remains the case of $a=1$ and $b$ composite, so $b \ge 9$.
We use \ref{psl5+uni} and \ref{psu5+uni} for the connected components of $F^\ast(M_2)$ and $F^\ast(M_3)$ and get $x$ into the big connected component.
\end{bew}

Before we can handle $B_n(q)$ we need $D_n(q)$ and ${}^2D_n(q)$.

\begin{lemma}
\label{pop8+uni}
Let $G \cong D_n(q)$ for $n \ge 4$. Then Theorem \ref{theorem_small_connected_components} holds.
\end{lemma}

\begin{bew}
If $q$ is odd, by \ref{odd_pcc} there exists a big connected component, containing all elements of order $p$.
There exists a subgroup $M_1$ of type $O_3(q) \perp O_{2n-3}(q)$, so if $x$ is not in the big connected component, $r \mid (q^n-1)(q^{n-1}-1)(q^{n-1}+1)$.\\
If $q$ is even, by Theorem \ref{theorem_connected_conjugacy_class}, there exists a connected component containing all elements of prime order $r$ for $r$ some divisor of
$q^2-1$. Let $M_1$ in class ${\cal C}_1$ of type $O^-_2(q) \perp O^-_{2n-2}(q)$.
By the structure of $M_1$, elements of order $s$ are in the big connected component, if $s$ is an odd prime divisor of $|\Omega^-_{2n-2}(q)|$.\\ 
So remain odd primes $r$, which divide $(q^n-1)(q^{n-1}-1)$.

Let $n$ even. 

Suppose $r \mid q^n-1$. 
If $q$ is odd, then $q^n-1 | |\Omega_{n+1}(q)| $ and $n+1 \le 2n-3$.
If $q$ is even, then $q^n-1 | |\Omega^-_{n+2}| $ and $n+2 \le 2n-2$.
This implies, that $x$ is in the big connected component by $M_1$ in both cases. 

Suppose $r \mid q^{n-1}-1$. A torus of type $q^{n-1}-1$ can be found in a subgroup $M_2$ of type $\GL_{n}(q).2$ in class ${\cal C}_2$. The structure of $M_2$ is described by Proposition 4.2.7 of \cite{KL}. If $q-1$ is not a 2-power, then $Z(F^\ast(M_2))$ contains elements of odd order, so $x$ is in the big connected component. 

We use \ref{psl4uni} and \ref{psl5+uni} for the connected components of $M_2$, if $q-1$ is a  2-power. 
Therefore $n-1$ is a prime. By observation of the action on the natural module we conclude, that this gives a small connected component.

Suppose $r \mid q^{n-1}+1$, so $q$ is odd. A torus of type $q^{n-1}+1$ is contained in a subgroup $M_3$ of
type $\GU_n(q)$ in class ${\cal C}_3$. The structure of $M_3$ is described by Proposition 4.3.18 of \cite{KL}. If $q+1$ is not a 2-power,
$Z(F^\ast(M_3))$ contains elements of odd order and $x$ is in the big connected component. We use \ref{psu4uni}
and \ref{psu5+uni} for the connected component of $M_3$, if $q+1$ is a 2-power. Therefore $n-1$ is a prime and 
we get again a small connected component.

Let $n$ odd.

Suppose $r \mid q^n-1$. A torus of type $q^{n}-1$ can be found in a subgroup $M_4$ of type $\GL_{n}(q).2$ in class ${\cal C}_2$. 
The structure of $M_4$ is described by Proposition 4.2.7 of \cite{KL}. 
If $q-1$ is not a 2-power, then $Z(F^\ast(M_4))$ contains elements of odd order, so $x$ is in the big connected component. 

We use \ref{psl4uni} and \ref{psl5+uni} for the connected components of $M_4$, if $q-1$ is a 2-power. 
Therefore $n$ is a prime and we get again a small connected component. 

Suppose $r \mid q^{n-1}-1$. 
If $q$ is odd, then $q^{n-1}-1 | |\Omega_{n}(q)| $ and $n \le 2n-3$.
If $q$ is even, then $q^{n-1}-1 | |\Omega^-_{n+1}| $ and $n+1 \le 2n-2$.
This implies, that $x$ is in the big connected component by $M_1$ in both cases.

Suppose $r \mid q^{n-1}+1$, so $q$ is odd and $n-1$ is even. Let $n-1=a \cdot b$ with $a>1$ a 2-power and $b$ odd. 
We can find a torus of type $q^{n-1}+1$ in a subgroup $M_5$ of type $\GU_b(q^a)$. 
We can find this subgroup in the following chain of subgroups: $\GU_b(q^a) \le O^-_{2b}(q^a) \le O^-_2(q) \perp O^-_{2n-2}(q) \le G$.
If $q+1$ is not a 2-power, then elements of order $r$ commute with elements of order $s$ for $s$ some divisor of $q+1$,
as visible in $O^-_2(q) \perp O^-_{2n-2}(q)$. 

If $q+1$ is a 2-power, then we have to analyze $M_5$. If $b=1$, we get a torus of type $q^{n-1}+1$ centralized by a 2-group of size $q+1$.
Observation of $x$ in the natural module shows, that we have a small connected component. 
If $b \ne 1$, the structure of $M_5$ as subgroup of $O^-_{2b}(q^a)$ is described by Proposition 4.3.18 of \cite{KL}. 
By Zsigmondy, $Z(F^\ast(M_5))$ contains elements of odd order $s$ with $s \mid q^a+1$. 
As $q^{2a}+1 \mid |\Omega_{2a+1}(q)|$, $a \le \frac{n}{3}$ and $n \ge 4$, we have $2a+1 \le 2n-3$, so $x$ is in the big connected
component by $M_1$.
\end{bew}

\begin{lemma}
\label{pom8+uni}
Let $G \cong {}^2D_n(q)$ for $n \ge 4$. Then Theorem \ref{theorem_small_connected_components} holds. 
\end{lemma}

\begin{bew}
If $q$ is odd, by \ref{odd_pcc} there exists a big connected component, containing all elements of order $p$.

There exists a subgroup $M_1$ of type $O_3(q) \perp O_{2n-3}(q)$, so if $x$ is not contained in the big connected component, 
$r \mid (q^n+1)(q^{n-1}-1)(q^{n-1}+1)$.

If $q$ is even, by Theorem \ref{theorem_connected_conjugacy_class}, there exists a connected component containing all elements of prime order $r$
for $r$ some prime divisor of $q^2-1$. 

Let $M_1$ in class ${\cal C}_1$ be of type $O^-_2(q) \perp O^+_{2n-2}(q)$. 
By the structure of $M_1$, elements of order $r$ are in that connected component, if $r \mid |\Omega^+_{2n}(q)|$,
so remain primes $r$, which divide $(q^n+1)(q^{n-1}+1)$.

Let $n$ even. \\
Suppose $r \mid q^n+1$. Let $n=a \cdot b$ with $a$ a 2-power and $b$ odd. Notice, $a \ne 1$ and $b=1$ gives 
a small connected component, as the torus of type $q^n+1$ is self centralizing. 

A torus of type $q^n+1$ is contained in a subgroup $M_2$ of type $\GU_b(q^a)$, which is contained
in a subgroup $M_3$ of type $O^-_{2b}(q^a)$. The structure of $M_2$ as subgroup of $M_3$ is described by Proposition 4.3.18 of \cite{KL}.
By Zsigmondy, $Z(F^\ast(M_2))$ contains elements of odd order $s$ with $s \mid q^a+1$. 
As $F^\ast(M_2)$ contains a $\PSL_2(q)$-section, $x$ is in the big connected component.

Suppose $r \mid q^{n-1}-1$, so $q$ is odd. A torus of type $q^{n-1}-1$ can be found in a subgroup $M_4$ of type 
$O^-_2(q) \perp O^+_{2n-2}(q)$. If $q+1$ is not a 2-power, $Z(F^\ast(M_4))$ containes elements of odd order,
so $x$ is in the big connected component. Else we may use \ref{pop8+uni} for the connected components of $F^\ast(M_4)$.
Observation of $x$ in the natural module shows, that we get a small connected component.

Suppose $r \mid q^{n-1}+1$. A torus of type $q^{n-1}+1$ can be found in a subgroup $M_5$ of type 
$O^+_2(q) \perp O^-_{2n-2}(q)$. If $q-1$ is not a 2-power, $Z(F^\ast(M_5))$ containes elements of odd order,
so $x$ is in the big connected component. Else we use induction for the connected components of $F^\ast(M_5)$. 
This gives another small connected component. 

Let $n$ odd. Suppose $r \mid q^n+1$. A torus of type $q^n+1$ can be found in a subgroup $M_6$ of type $\GU_n(q)$.
The structure of $M_6$ is described by Proposition 4.3.18 of \cite{KL}. If $q+1$ is not a 2-power,
then $Z(F^\ast(M_6))$ contains elements of odd order and $x$ is contained in the big connected component.
We use \ref{psu4uni} and \ref{psu5+uni} for the connected components of $F^\ast(M_6)$. This gives a small connected component. 

Suppose $r \mid q^{n-1}-1$, so $q$ is odd. As $q^{n-1}-1 \mid |\Omega_n(q)|$ and $n \le 2n-3$, $x$ is in the big connected component by $M_1$.

Suppose $r \mid q^{n-1}+1$. A torus of type $q^{n-1}-1$ can be found in a subgroup $M_7$ of type $O_2^+(q) \perp O^-_{2n-2}(q)$.
If $q-1$ is not a 2-power, then $Z(F^\ast(M_7))$ contains elements of odd order, so $x$ is in the big connected component.

Else we get the structure of the connected components of $F^\ast(M_7)$ by induction. This gives another small connected component.
\end{bew}

\begin{lemma}
\label{po7+uni}
Let $G \cong B_n(q)$ for $n \ge 3$, so $q$ is odd. Then Theorem \ref{theorem_small_connected_components} holds. 
\end{lemma}

\begin{bew}
By \ref{odd_pcc} there exists a big connected component, containing all elements of order $p$.
There exist subgroups $M_1$ of type $O_3(q) \perp O^+_{2n-2}(q)$ and $M_2$ of type $O_3(q) \perp O^-_{2n-2}(q)$ 
so either $x$ is in the big connected component or $r \mid (q^n-1)(q^n+1)$.

Suppose $r \mid q^n-1$. A torus of type $q^n-1$ can be found in a subgroup $M_3$ of type $O_1(q) \perp O^+_{2n}(q)$.
We use \ref{pop8+uni}  for the structure of the connected components of $F^\ast(M_3)$, but restrict to cases not contained in $M_1$.
Observation of the natural module gives a small connected component. 

Suppose $r \mid q^n+1$. A torus of type $q^n+1$ can be found in a subgroup $M_4$ of type $O_1(q) \perp O^-_{2n}(q)$.
We use \ref{pom8+uni}  for the structure of the connected components of $F^\ast(M_4)$, but restrict to cases bot contained in $M_1$. 
Observation of the action of corresponding elements on the natural module shows, that we get two more small connected components.
\end{bew}

\begin{lemma}
\label{G2uni}
Let $G \cong G_2(q)$. Then Theorem \ref{theorem_small_connected_components} holds. 
\end{lemma}

\begin{bew}
By Theorem \ref{theorem_big_component}, $q>2$. 

If $q$ odd, by \ref{odd_pcc} there exists a big connected component, containing all elements of order $p$.\\
If $q$ even, by Theorem \ref{theorem_connected_conjugacy_class}, there is a connected conjugacy class $y^G$, $o(y)=r$ for $r$ some prime divisor of $q^2-1$ ($r=3$ for $q=4$). 
By \cite{LSS} there exists a subgroup $M_1$ of type $\SL_2(q) \circ \SL_2(q)$. Therefore $d_q(r) \in \{3,6\}$. 

Suppose $d_q(r)=3$. By \cite{LSS} there exists a subgroup $M_2$ of type $\SL_3(q)$, which has a nontrivial center, if $3 \mid q-1$.
Suppose $d_q(r)=6$. By \cite{LSS} there exists a subgroup $M_3$ of type $\SU_3(q)$, which has a nontrivial center, if $3 \mid q+1$.
By \cite{CLSS} we get a list of maximal subgroups, which contain all centralizers of elements in $G$. Using this list we conclude,
that these conditions indeed give small connected components.
\end{bew}

\begin{lemma}
\label{3D4uni}
Let $G \cong {}^3D_4(q)$. Then Theorem \ref{theorem_small_connected_components} holds. 
\end{lemma}

\begin{bew}
If $q$ is odd, by \ref{odd_pcc} there exists a big connected component, containing all elements of order $p$.\\
If $q$ is even, by Theorem \ref{theorem_connected_conjugacy_class}, there is a connected conjugacy class $y^G$, $o(y)=s$ for $s$ some prime divisor of $q^2-1$.
By \cite{LSS} there exists a subgroup $M_1$ of type $\SL_2(q) \circ \SL_2(q^3)$. Therefore $d_q(r)=12$. As the torus of size $q^4-q^2+1$ is self centralizing
by \cite{K3D4}, we get a small connected component. 
\end{bew}

\begin{lemma}
\label{2F4uni}
Let $G \cong {}^2F_4(q)'$ for $q$ even. Then Theorem \ref{theorem_small_connected_components} holds.
\end{lemma} 

\begin{bew}
For $q=2$ we use the centralizer sizes in \cite{ATLAS}. 

Recall, that $3 \mid q+1$ and $5 \mid q^2+1$, as $q$ is an odd power of $2$. 

By Theorem \ref{theorem_connected_conjugacy_class}, there is a connected conjugacy class $y^G$, $o(y)=r$ for $r$ some prime divisor of $q^2+1$.  
By \cite{Malle}, subgroups of type $\SU_3(q)$, ${}^2B_2(q) \times {}^2B_2(q)$ and $\Sp_4(q)\ge \PSL_2(q) \times \PSL_2(q)$ exist.
Let $\rho$ be the set of element orders of ${\cal C}_y$. 
We have $\pi(q^2+1) \subseteq \rho$, so $\pi(q-1) \subseteq \rho$ from ${}^2B_2(q) \times {}^2B_2(q)$,
so $\pi(q+1) \subseteq \rho$ from $\PSL_2(q) \times \PSL_2(q)$, so $3 \in \rho$, so $\pi(q^3+1) \subseteq \rho$ from $\SU_3(q)$. 

As self centralizing subgroups of size $q^2+\sqrt{2q^3}+q+\sqrt{2q}+1$ and $q^2-\sqrt{2q^3}+q-\sqrt{2q}+1$ exist
with $(q^2+\sqrt{2q^3}+q+\sqrt{2q}+1)(q^2-\sqrt{2q^3}+q-\sqrt{2q}+1)=q^4-q^2+1= \Phi_{12}(q)$, these subgroups produce small connected components. 
\end{bew}

\begin{lemma}
\label{F4uni}
Let $G \cong F_4(q)$. Then Theorem \ref{theorem_small_connected_components} holds. 
\end{lemma}

\begin{bew}
If $q$ is odd, by \ref{odd_pcc} there exists a big connected component, containing all elements of order $p$.

If $q$ is even, by Theorem \ref{theorem_connected_conjugacy_class}, we have a connected component containing all elements of order $r$ for $r$ some divisor
of $q^2-1$. \\
By \cite{LSS} there exist subgroup $M_1$ of type $\Omega_9(q)$ ($q$ odd) or $\Sp_8(q)$ ($q$ even) and $M_2$ of type ${}^3D_4(q)$.
From the group order formula, $r \mid |M_1||M_2|$. By \ref{po7+uni},\ref{psp6+uni} and \ref{3D4uni}, $d_q(r) \in \{ 8,12 \}$.
By \cite{CLSS} we get a list of maximal subgroups, which contain all centralizers of elements in $G$. Using this list we conclude,
that these conditions indeed give small connected components.
\end{bew}

\begin{lemma}
\label{E678uni}
Let $G \cong E_6(q)$, ${}^2E_6(q)$, $E_7(q)$ or $E_8(q)$. Then Theorem \ref{theorem_small_connected_components} holds. 
\end{lemma}

\begin{bew}
We first reduce to a few cases and later use the list of maximal subgroups given in \cite{CLSS},
which contain centralizers of all elements in $G$. Using this list we conclude, that our conditions indeed give small connected components.\\

If $q$ is odd, by \ref{odd_pcc} there exists a big connected component, containing all elements of order $p$.\\
If $q$ is even, by Theorem \ref{theorem_connected_conjugacy_class}, there exists a connected component containing all elements of order $r$ for
$r$ some prime divisor of $q^2-1$.\\
Consider $G \cong E_6(q)$. By \cite{LSS} there exists a subgroup of type $\SL_2(q) \circ \SL_6(q)$.
Therefore $d_q(r) \in \{8,9,12\}$. 

Suppose $d_q(r)=12$. By \cite{LSS} there exists a subgroup of type $({}^3D_4(q) \circ \frac{q^2+q+1}{(q-1,3)}$,
therefore $x$ is in the big connected subgroup. 

Suppose $d_q(r)=8$. By \cite{LSS} there exists a subgroup of type $\Omega_{10}^+(q) \times \frac{(q-1)}{(q-1,3)}$. 
We use \ref{pop8+uni} for the connected components of this group, if $\frac{q-1}{(q-1,3)}$ is a 2-power. 
We get $x$ centralized by a subgroup of size $(q+1) \frac{q-1}{(q-1,3)}$. This is a 2-power, iff $q=3$ or $q=7$ by \ref{qLemma}.

The case $d_q(r)=9$ gives a small connected component.\\
Consider $G \cong {}^2E_6(q)$. By \cite{LSS} there exists a subgroup of type $\SL_2(q) \circ \SU_6(q)$.
Therefore $d_q(r) \in \{8,12,18\}$. 

Suppose $d_q(r)=12$. By \cite{LSS} there exists a subgroup of type $({}^3D_4(q) \circ \frac{q^2-q+1}{q+1,3})$,
therefore $x$ is in the big connected subgroup, except $q=2$. 

Suppose $d_q(r)=8$. By \cite{LSS} there exists a subgroup of type $\Omega_{10}^-(q) \times \frac{q+1}{(q+1,3)}$. 
We use \ref{pom8+uni} for the connected components of this group, if $\frac{q+1}{(q+1,3)}$ is a 2-power.
We get $x$ centralized by a subgroup of size $(q-1) \frac{q+1}{(q+1,3)}$. This is a 2-power, iff $q=2,3$ or $q=5$ by \ref{qLemma}.
 
The case $d_q(r)=18$ gives a small connected component.\\

Consider $G \cong E_7(q)$. By \cite{LSS} there exists a subgroup of type $\SL_2(q) \circ \Omega^+_{12}(q)$. 
This gives $d_q(r) \in \{7,9,12,14,18\}$. By \cite{LSS} there exists a subgroup of type $\PSL_2(q^3) \times {}^3D_4(q)$,
which shows $x$ in the big connected component for $d_q(r)=12$. \\
A subgroup of type $\PSL_2(q^7)$ give small connected for $d_q(r) \in \{ 7,14\}$, if $q-1$ resp. $q+1$ is a 2-power.\\
Subgroups of type $E_6(q) \circ (q-1)$ and ${}^2E_6(q) \circ (q+1)$ give more small connected components for $d_q(r) \in \{ 9,18 \}$ and $q-1$ resp. $q+1$ a 2-power. 
For the existence of these subgroups we use \cite{LSS}. \\

Consider $G \cong  E_8(q)$. By \cite{LSS} there exists a subgroup of type $\SL_2(q) \circ E_7(q)$.
This gives $d_q(r) \in \{ 15,20,24,30\}$. By \cite{LSS} there exists
a subgroup of type $\SU_5(q^2)$, which contains a torus of type $\frac{q^{10}+1}{q^2+1}$ and 
has a nontrivial center, if $5 \mid q^2+1$. This completes the list of small connected components.
\end{bew}


\begin{thebibliography}{99}

\bibitem[AAM]{AAM} A. Abdollahi, S. Akbari, H.R. Maimani, Non-commuting graph of a group, {\em J. Algebra} {\bf 298}
(2006), 468-492.

\bibitem[AS]{AS} M. Aschbacher and Y. Segev, The uniqueness of groups of Lyons type, {\em J. Amer. Math. Soc.} {\bf 5} (1992), 75-98.

\bibitem[BBPR]{BBPR} C. Bates, D. Bundy, S. Perkins, P. Rowley, Commuting involution graphs for symmetric groups,
{\em J. Algebra} {\bf 266} no. 1 (2003), 133-153.

\bibitem[BS]{BS} B. Baumeister, A.Stein,  The finite Bruck loops, Preprint 2009.

\bibitem[B]{bender} H. Bender, Transitive Gruppen gerader Ordnung, in denen jede Involution genau einen Punkt festl\"a{\ss}t, {\em J. Algebra} {\bf 17} (1971), 527-54.

\bibitem[Car]{Car} R. Carter, {Simple Groups of Lie Type}, John Wiley and Sons, London, 1972. 

\bibitem[Car2]{Car2} R. Carter, {Finite Groups of Lie Type: conjugacy classes and complex characters}, Wiley-Interscience, 1985. 

\bibitem[Coo]{Coo} B.N. Cooperstein, The maximal subgroups of $G_2(2^n)$. {\em J. Algebra} {\bf 70}, 23-36 (1981).

\bibitem[CLSS]{CLSS} A.M. Cohen, M.W. Liebeck, J. Saxl, G.M. Seitz, The local maximal subgroups of exceptional groups of Lie type, finite and algebraic
{\em Proc. Lond. Math. Soc.}, III. Ser. {\bf 64}, No.1 (1992), 21-48.
  
\bibitem[ATLAS]{ATLAS} J.H. Conway, R.T. Curtis, S.P. Norton, R.A. Parker ,R.A. Wilson {An ATLAS of finite groups} Oxford University Press, 1985.

\bibitem[Fi]{Fi} B. Fischer, Finite groups generated by $3$-transpositions, I, {\em Invent. Math}, {\bf 13} (1971), 232-246.

\bibitem[GLS3]{GLS3}  D. Gorenstein, R. Lyons, R. Solomon {The Classification of Finite Simple Groups, Volume 3}, Math. Surveys and Monographs, Vol. 40, No. 3, American Math. Soc., 1998.

\bibitem[IJ1]{IJ1} A. Iranmanesh, A. Jafarzadeh, Characterization of finite groups by their commuting graph,
{\em Acta Math. Acad. Paedagog. Nyhazi} {\bf 23} no. 1 (2007), 7-13.

\bibitem[IJ2]{IJ2} A. Iranmanesh, A. Jafarzadeh, On the commuting graph associated with the symmetric an alternating
groups, {\em J. Algebra Appl.} {\bf 7} no. 1 (2008), 129-146.

\bibitem[K2G2]{K2G2} P. Kleidman The maximal subgroups of the Chevalley groups $G_2(q)$with $q$ odd, of the Ree groups ${}^2G_2(q)$ and of their automorphism groups
{\em J.Algebra} {\bf 117} (1988), 30-71.

\bibitem[K3D4]{K3D4} P. Kleidman {The maximal subgroups of the Steinberg triality groups ${}^3D_4(q)$ and of their automorphism groups} {\em J.Algebra} {\bf 115} (1988), 182-199.

\bibitem[KL]{KL} P. Kleidman, M. Liebeck {The Subgroup Structure of the Finite Classical Groups} LMS Lecture Note Series 129, Cambridge University Press, 1990.

\bibitem[LSS]{LSS} M.W. Liebeck, J. Saxl, G.M. Seitz, {Subgroups of maximal rank in finite exceptional groups of Lie type} 
{\em Proc. Lond. Math. Soc.}, III. Ser. {\bf 65}, No.2 (1992), 297-325.

\bibitem[Malle]{Malle} G. Malle, {The maximal subgroups of ${}^2F_4(q^2)$}, {\em J.Algebra}  {\bf 139} (1991), 52-69.

\bibitem[Pe]{Pe} S. Perkins, {Commuting involution graphs for $\tilde{A}_n$}, {\em Arch. Math.} {\bf 86} no. 1 (2006), 16-25.

\bibitem[S]{S} A. Stein, {On Bruck Loops of 2-power Exponent, II}, Preprint 2009.

\bibitem[Se]{Se} Y. Segev, {On finite homomorphic images of the multiplicative group of a division algebra.} {\em Ann. of Math.} (2) {\bf 149} , no. 1 (1999), 219-251.

\bibitem[St]{Stellmacher} B. Stellmacher {Einfache Gruppen, die von einer Konjugiertenklasse von Elementen der Ordnung 3 erzeugt werden.}
{\em J.Algebra}  {\bf 30} (1974), 320-354.

\bibitem[Sz1]{suz-CA-groups} M. Suzuki, {The nonexistence of a certain type of simple groups of odd order}, {\em Proc. AMS} {\bf 8} (1957), 686-695

\bibitem[Sz]{suzuki} M. Suzuki, {Group Theory} Vol.I,II , Springer 1982,1986

\bibitem[W]{weisner} L. Weisner {Groups, in which the Normaliser of every element except the identity is abelian}, {\em Bull. AMS} {\bf 31} (1925), 413-416

 
\end{thebibliography}
\end{document}